\RequirePackage{ifpdf}
\ifpdf % We are running pdfTeX in pdf mode
\documentclass[pdftex]{sigma}
\else
\documentclass{sigma}
\fi

\numberwithin{equation}{section}

\begin{document}

%\allowdisplaybreaks

\renewcommand{\PaperNumber}{076}

\FirstPageHeading

\ShortArticleName{The Symmetrical $H_{q}$-Semiclassical Orthogonal
Polynomials of Class One}

\ArticleName{The Symmetrical $\boldsymbol{H_{q}}$-Semiclassical\\ Orthogonal
Polynomials of Class One}

\Author{Abdallah GHRESSI~$^\dag$ and Lotf\/i KH\'{E}RIJI~$^\ddag$}

\AuthorNameForHeading{A. Ghressi and L. Kh\'{e}riji}

\Address{$^\dag$~Facult\'{e} des Sciences de Gab\`{e}s,
Route de Mednine 6029 Gab\`{e}s, Tunisia}
\EmailD{\href{mailto:Abdallah.Ghrissi@fsg.rnu.tn}{Abdallah.Ghrissi@fsg.rnu.tn}}

\Address{$^\ddag$~Institut Sup\'{e}rieur des Sciences Appliqu\'{e}es
et de Technologies de Gab\`{e}s,\\
\hphantom{$^\ddag$}~Rue Omar Ibn El-Khattab 6072 Gab\`{e}s, Tunisia}
\EmailD{\href{mailto:Lotfi.Kheriji@issatgb.rnu.tn}{Lotfi.Kheriji@issatgb.rnu.tn}}

\ArticleDates{Received December 12, 2008, in f\/inal form July 07, 2009;  Published online July 22, 2009}

\Abstract{We investigate the quadratic decomposition and duality to
classify symmetrical $H_{q}$-semiclassical orthogonal
$q$-polynomials of class one where $H_{q}$ is the Hahn's operator.
For any canonical situation, the recurrence coef\/f\/icients, the
$q$-analog of the distributional equation of Pearson type, the
moments and integral or discrete representations are given.}


\Classification{33C45; 42C05}

\section{Introduction}
Orthogonal polynomials (OP) have been a subject of research in the
last hundred and f\/ifty years. The orthogonality considered in our
contribution is related to a form (regular linear functional)~\cite{8,24}
and not only to a positive measure. By classical orthogonal
polynomials sequences (OPS), we refer to Hermite, Laguerre, Bessel
and Jacobi polynomials. In the literature, the extension of
classical (OPS) can be done from dif\/ferent approaches such that the
hypergeometric character \cite{7,8,11,18,22} and the distributional
equation of Pearson type \cite{6,8,20,29,32}. A natural generalization of
the classical character is the semiclassical one introduced by J.A.~Shohat in~\cite{35}. This theory was developed by P.~Maroni and
extensively studied by P.~Maroni and coworkers in the last decade~\cite{1,24,26,28,32}.
Let $\Phi$ monic and $\Psi$ be two polynomials,
$\deg \Phi=t\geq 0$,  $\deg \Psi=p\geq 1$. We suppose that the pair
($\Phi,\Psi$) is admissible, i.e., when $p=t-1$, writing
$\Psi(x)=a_{p}x^{p}+\cdots$, then $a_{p}\neq n+1$, $n\in \mathbb{N}$.
A form $u$ is called semiclassical when it is regular and satisf\/ies
the distributional equation of Pearson type
\begin{gather}
D(\Phi u) +\Psi u=0, \label{1.1}
\end{gather}
where the pair($\Phi,\Psi$) is admissible and $D$ is the derivative
operator. The corresponding monic orthogonal polynomials sequence
(MOPS) $\{B_{n}\}_{n\geq 0}$ is called semiclassical. Moreover, if
$u$ is semiclassical satisfying \eqref{1.1}, the class of $u$, denoted $s$
is def\/ined by
\begin{gather*}
s:=\min \big(\max(\deg \Phi-2,\deg \Psi-1)\big)\geq 0, %\label{1.2}
\end{gather*}
where the minimum is taken over all pairs ($\Phi,\Psi$) satisfying~\eqref{1.1}. In particular, when $s=0$ the classical case is recovered.

Symmetrical semiclassical forms of class one are well
described in~\cite{1}, see also~\cite{6}; there are three canonical
situations:
\begin{enumerate}\itemsep=0pt
\item[1)] The generalized Hermite form $\mathcal{H}(\mu)$
($\mu\neq0$, $\mu\neq-n-\frac{1}{2}$,
$n\geq0$) satisfying the distributional equation of
Pearson type
\begin{gather}
D(x\mathcal{H}(\mu))+\bigl(2x^{2}-(2\mu+1)\bigr)\mathcal{H}(\mu)=0.\label{1.3}
\end{gather}

\item[2)] The generalized Gegenbauer
$\mathcal{G}(\alpha,\beta)$
($\alpha\neq-n-1$, $\beta\neq-n-1$,
$\beta\neq-\frac{1}{2}$,
$\alpha+\beta\neq-n-1$, $n\geq0$)
satisfying the distributional equation of Pearson type
\begin{gather}
D\big(x\big(x^{2}-1\big)\mathcal{G}(\alpha,\beta)\big)+
\big(-2(\alpha+\beta+2)x^{2}+2(\beta+1)\big)\mathcal{G}(\alpha,\beta)=0.
\label{1.4}
\end{gather}
For further properties of the generalized Hermite and the
generalized Gegenbauer polynomials see
\cite{8,15,26}.
\item[3)] The form $\mathcal{B}[\nu]$ of Bessel kind
($\nu\neq-n-1$, $n\geq0$) \cite{1,31} satisfying
the distributional equation of Pearson type
\begin{gather}
D\big(x^{3}\mathcal{B}[\nu]\big)-\big(2(\nu+1)x^{2}+\tfrac{1}{2}\big)\mathcal{B}[\nu]=0.\label{1.5}
\end{gather}
For an integral representation of $\mathcal{B}[\nu]$ and some
additional features of the associated (MOPS) see~\cite{14}.
\end{enumerate}

Other families of semiclassical orthogonal polynomials of class
greater than one were discovered by solving functional equations of
the type $P (x) u=Q (x)v$, where $ P$, $Q$ are two polynomials
cunningly chosen and $u$, $v$ two linear forms \cite{19,26,27,34}. For
other relevant works
in the semiclassical case see \cite{5,23}.

In \cite{21}, instead of the derivative operator, the $q$-dif\/ference one
is used to establish the theo\-ry and characterizations of
$H_{q}$-semiclassical orthogonal $q$-polynomials. Some examples of
$H_{q}$-semiclassical orthogonal $q$-polynomials are given in
\cite{2,13}. The $H_{q}$-classical case is exhaustively described in
\cite{20,32}. Moreover, in \cite{30} the symmetrical
$D_{\omega}$-semiclassical orthogonal polynomials of class one are
completely described by solving the system of their Laguerre--Freud
equations where $D_{w}$ is the Hahn's operator.

So, the aim of this paper is to present the classif\/ication of the
symmetrical $H_{q}$-semiclassical orthogonal $q$-polynomials of
class one by investigating the quadratic operator $\sigma$, the
$q$-analog of the distributional equation of Pearson type satisf\/ied
by the corresponding form and some $H_{q}$-classical situations (see
Tables~\ref{table1} and~\ref{table2}) in connection
with our problem. Among the obtained canonical cases, three are well
known: two symmetrical Brenke type (MOPS)~\cite{8,9,10} and a~symmetrical
case of the Al-Salam and Verma (MOPS)~\cite{2}. Also, $q$-analogues of
$\mathcal{H}(\mu)$, $\mathcal{G}(\alpha,\beta)$ and
$\mathcal{B}[\nu]$ appear. In \cite{3,33}, the authors have established,
up a dilation, a $q$-analogues of $\mathcal{H}(\mu)$ and
$\mathcal{B}[\nu]$ using other methods. For any canonical case, we
determine the recurrence coef\/f\/icient, the $q$-analog of the
distributional equation of Pearson type, the moments and a discrete
measure or an integral representation.

\begin{table}[t]\centering\small
\caption{Canonical cases.} \label{table1}

\vspace{1mm}

\begin{tabular}{|p{15.5cm}|}
  \hline
   \tsep{1mm}\bsep{1mm} $H_{q}$-classical linear form \\
   \hline
\tsep{1mm}\bsep{1mm} case 1.1\hspace{1cm}  ${\mathcal{U}}$ \\
   \hline
\tsep{2mm} $\widehat{\beta}_{n}=\{1-(1+q)q^{n}\}q^{n-1}$, $n\geq 0$, \\
\tsep{1mm} $\widehat{\gamma}_{n+1}=(q^{n+1}-1)q^{3n}$, $n\geq 0$, \\
\tsep{1mm} $H_{q}(x\mathcal{U})-(q-1)^{-1}(x+1) \mathcal{U}=0$,\\
\tsep{1mm} $(\mathcal{U})_{n}=(-1)^{n}q^{\frac{1}{2}n(n-1)}$, $n\geq 0$,\\
\tsep{1mm} $\mathcal{U}= \displaystyle\sum_{k=0}^{\infty} (-1)^{k}
\frac{q^{-k^{2}} s(k)}{(q^{-1};q^{-1})_{k}} \delta_{-q^{k}}$, $q>1$,\\
\tsep{1mm} where $s(k)=\displaystyle\sum_{m=0}^{\infty}\frac{q^{-(\frac{1}{2}m(m+1)+km)}}
{(q^{-1};q^{-1})_{k}}\varepsilon_{m+k}$,
$\varepsilon_{2k}=(q-1)^{k}$, $k\geq0$  and
$\varepsilon_{2k+1}=0$, $k\geq0$ \bsep{1mm}\\
  \hline
\tsep{1mm} case 1.2\hspace{1cm}  little $q$-Laguerre ${\mathcal{L}(a,q)}$
\big($a\neq 0$, $a\neq q^{-n-1}$, $n\geq0$\big) \bsep{1mm}\\
\hline
\tsep{2mm}
$\widehat{\beta}_{n}=\{1+a-a(1+q)q^{n}\}q^{n}$, $n\geq 0,$ \\
\tsep{1mm} $\widehat{\gamma}_{n+1}=a(1-q^{n+1})(1-aq^{n+1})q^{2n+1}$, $n\geq 0$, \\
$H_{q}(x\mathcal{L}(a,q))-(aq)^{-1}(q-1)^{-1}\{x-1+aq\} \mathcal{L}(a,q)=0$,
\\
\tsep{1mm}  $(\mathcal{L}(a,q))_{n}=(aq;q)_{n}$, $n \geq 0$,\\
\tsep{1mm}  $ \mathcal{L}(a,q)=(aq;q)_{\infty}\displaystyle\sum_{k=0}^{\infty}
  \displaystyle\frac{(aq)^{k}}
  {(q;q)_{k}}\delta_{q^{k}}$, $0<q<1$, $0<a<q^{-1}$,\\
\tsep{1mm} $ \langle
\mathcal{L}(a,q),f\rangle=K\displaystyle\int_{0}^{q^{-1}}x^{\frac{\ln
a}{\ln
q}}(qx;q)_{\infty}f(x)dx$, $f\in \mathcal{P}$, $0<q<1$, $0<a<q^{-1}$,\\
\tsep{1mm} where  $K^{-1}=q^{-\frac{\ln a}{\ln
q}-1}\displaystyle\int_{0}^{1}x^{\frac{\ln a}{\ln
q}}(x;q)_{\infty}dx$,\\
\tsep{1mm} $\mathcal{L}(a,q)=\frac{1}{(a;q^{-1})_{\infty}}
\displaystyle\sum_{k=0}^{\infty}
\frac{q^{-\frac{1}{2}k(k-1)}}{(q^{-1};q^{-1})_{k}}(-a)^{k}
\delta_{q^{k}}$, $q>1$, $a<0$\bsep{1mm}\\
 \hline
\tsep{1mm}  case 1.3\hspace{1cm}  Wall  ${\mathcal{W}(b,q)}$
\big($b\neq 0$, $b\neq q^{-n}$, $n\geq 0$\big)\bsep{1mm}
 \\
 \hline
\tsep{2mm} $\widehat{\beta}_{n}=\{b+q-b(1+q)q^{n}\}q^{n}$, $n\geq 0$, \\
\tsep{1mm} $\widehat{\gamma}_{n+1}=b(1-q^{n+1})(1-bq^{n})q^{2n+2}$, $n\geq 0,$ \\
\tsep{1mm} $H_{q}(x\mathcal{W}(b,q))-b^{-1}(q-1)^{-1}(q^{-1}x+b-1)\mathcal{W}(b,q)=0$,\\
\tsep{1mm} $(\mathcal{W}(b,q))_{n}=q^{n}(b;q)_{n}$, $n \geq 0$,\\
\tsep{1mm} $\langle\mathcal{W}(b,q),f\rangle=\frac{(b;q)_{\infty}}{2}
\displaystyle\sum_{k=0}^{\infty}\frac{b^{k}}{(q;q)_{k}} \langle
\delta_{q^{1+k}},f\rangle+\frac{K}{2}
\displaystyle\int_{0}^{1}x^{\frac{\ln b}{\ln q}-1}
(x;q)_{\infty}f(x)dx$, \\
 \tsep{1mm} where  $K^{-1}=\displaystyle\int_{0}^{1}x^{\frac{\ln b}{\ln
q}-1} (x;q)_{\infty}dx$, $f\in \mathcal{P}$, $0<q<1$, $0<b<1$,\\
$\langle\mathcal{W}(b,q),f\rangle=\frac{1}{(bq^{-1};q^{-1})_{\infty}}
\displaystyle\sum_{k=0}^{\infty}
\frac{q^{-\frac{1}{2}k(k+1)}}{(q^{-1};q^{-1})_{k}}(-b)^{k}\langle
\delta_{q^{1+k}},f\rangle$, $f\in \mathcal{P}$, $q>1$, $b \neq
q^{\pm k}$, $k\geq0$ \bsep{1mm}\\
\hline
\end{tabular}
\end{table}

\begin{table}[t]\centering\small
\centerline{\it Continuation of Table~{\rm \ref{table1}}.}\vspace{1mm}
 \begin{tabular}{|p{15.5cm}|}
  \hline
 \tsep{1mm}\bsep{1mm} case 1.4 \hspace{1cm} Generalized
$q^{-1}$-Laguerre ${\mathcal{U}^{(\alpha)}(b,q)}$
\big($b\neq 0$, $b\neq q^{n+1+\alpha}$, $n\geq 0$\big) \\
 \hline
\tsep{2mm}
$\widehat{\beta}_{n}=\{1-q^{-n-1}+q^{-1}(1-bq^{-n-\alpha})\}q^{2n+\alpha+1}$, $n\geq 0$, \\
\tsep{1mm} $\widehat{\gamma}_{n+1}=(1-q^{-n-1})(1-bq^{-n-1-\alpha})q^{4n+2\alpha+3}$, $n\geq 0$, \\
\tsep{1mm} $H_{q}(x\mathcal{U}^{(\alpha)}(b,q))+(q-1)^{-1}q^{-\alpha-1}(x+b-q^{\alpha+1})
\mathcal{U}^{(\alpha)}(b,q)=0$,
\\
\tsep{1mm} $(\mathcal{U}^{(\alpha)}(b,q))_{n}=(-b)^{n}(b^{-1}q^{\alpha+1};q)_{n}$,
$n \geq 0$,\\
\tsep{1mm} $\langle
\mathcal{U}^{(\alpha)}(b,q),f\rangle=(b^{-1}q^{\alpha+1};q)_{\infty}
\displaystyle\sum_{k=0}^{\infty}
\frac{(b^{-1}q^{\alpha+1})^{k}}{(q;q)_{k}}\langle
\delta_{-bq^{k}},f\rangle$,\\
\tsep{1mm} \hspace*{2.8cm} $f\in \mathcal{P}$,
$0<q<1$, $b >q^{\alpha+1}$, $\alpha \in\mathbb{R}$,\\
\tsep{1mm} $\langle\mathcal{U}^{(\alpha)}(b,q),f\rangle=
\frac{1}{2(b^{-1}q^{\alpha};q^{-1})_{\infty}}
\displaystyle\sum_{k=0}^{\infty}
\frac{q^{-\frac{1}{2}k(k-1)}(-b^{-1}q^{\alpha})^{k}}{(q^{-1};q^{-1})_{k}}
\langle \delta_{-bq^{k}},f\rangle$\\
$\displaystyle \phantom{\langle\mathcal{U}^{(\alpha)}(b,q),f\rangle=}{}
+\frac{K}{2}
\displaystyle\int_{0}^{\infty}\frac{x^{\alpha-\frac{\ln b}{\ln q}}}
{(-b^{-1}x;q^{-1})_{\infty}}f(x)dx$,
$f\in \mathcal{P}$, $q>1$, $q^{\alpha}<b<q^{\alpha+1}$, $\alpha \in\mathbb{R}$, \\
\tsep{1mm} where
$K^{-1}=\displaystyle\int_{0}^{\infty}\frac{x^{\alpha-\frac{\ln
b}{\ln q}}} {(-b^{-1}x;q^{-1})_{\infty}}dx$  is given by \eqref{2.10}\bsep{3mm} \\
\hline
\tsep{1mm}\bsep{1mm} case 1.5\hspace{1cm}   Alternative $q$-Charlier ${\mathcal{A}(a,q)}$  \big($a\neq 0$, $a\neq
-q^{-n}$, $n\geq0$\big)\\
\hline
\tsep{4mm}
$\displaystyle \widehat{\beta}_{n}=\frac{1+aq^{n-1}+aq^{n}-aq^{2n}}
{(1+aq^{2n-1})(1+aq^{2n+1})} q^{n}$, $n\geq 0,$ \\
\tsep{1mm} $\displaystyle \widehat{\gamma}_{n+1}=aq^{3n+1}\frac{(1-q^{n+1})(1+aq^{n})}{(1+aq^{2n})
(1+aq^{2n+1})^{2}(1+aq^{2n+2})}$, $n\geq 0,$ \\
\tsep{1mm} $H_{q}(x^{2}\mathcal{A}(a,q))-(aq)^{-1}(q-1)^{-1}\{(1+aq)x-1\}
\mathcal{A}(a,q)=0$,
\\
\tsep{1mm}  $\displaystyle (\mathcal{A}(a,q))_{n}=\frac{1}{(-aq;q)_{n}}$, $ n \geq 0$,\\
\tsep{1mm} $\displaystyle \langle
\mathcal{A}(a,q),f\rangle=\frac{1}{2(-aq;q)_{\infty}}\sum_{k=0}^{\infty}
  \displaystyle\frac{q^{\frac{1}{2}k(k+1)}a^{k}}
  {(q;q)_{k}}\langle\delta_{q^{k}},f\rangle +q^{\frac{1}{2}(\frac{\ln a}{\ln
q}+\frac{1}{2})^{2}}\frac{(-a^{-1};q)_{\infty}}{2 \sqrt{2\pi \ln
q^{-1}}}$\\
$\displaystyle \phantom{\langle \mathcal{A}(a,q),f\rangle=}{} \times \int_{0}^{\infty}x^{\frac{\ln a}{\ln
q}-\frac{1}{2}}(qx;q)_{\infty}\exp\left(-\frac{\ln^{2}x}{2\ln
q^{-1}}\right)
f(x)dx $,
$f\in \mathcal{P}$, $0<q<1$, $a > 0$\bsep{3mm}\\
\hline
\end{tabular}
\end{table}

\begin{table}[t]\centering\small
\centerline{\it Continuation of Table~{\rm \ref{table1}}.}\vspace{1mm}
 \begin{tabular}{|p{15.5cm}|}
\hline
\tsep{1mm}\bsep{1mm} case 1.6 \hspace{1cm}  little  $q$-Jacobi
${\mathcal{U}(a,b,q)}$  \big($ab\neq 0$, $a\neq q^{-n-1}$, $b\neq q^{-n-1}$, $ab\neq
q^{-n}$, $n\geq0$\big)\\
\hline
\tsep{4mm}
 $\displaystyle \widehat{\beta}_{n}=\frac{(1+a)(1+abq^{2n+1})-a(1+b)(1+q)q^{n}}
{(1-abq^{2n})(1-abq^{2n+2})} q^{n}$, $n\geq 0$, \\
\tsep{1mm}  $\displaystyle \widehat{\gamma}_{n+1}=aq^{2n+1}
\frac{(1-q^{n+1})(1-aq^{n+1})(1-bq^{n+1})(1-abq^{n+1})}
{(1-abq^{2n+1}) (1-abq^{2n+2})^{2}(1-abq^{2n+3})}$, $n\geq 0$, \\
\tsep{1mm}  $H_{q}(x(x-b^{-1}q^{-1})\mathcal{U}(a,b,q))
 +(abq^{2}(q-1))^{-1}\{(1-abq^{2})x+aq-1\} \mathcal{U}(a,b,q)=0$,
\\
$\displaystyle (\mathcal{U}(a,b,q))_{n}=\frac{(aq;q)_{n}}{(abq^{2};q)_{n}}$,
$n \geq 0$, \\
\tsep{1mm} $\langle
\mathcal{U}(a,b,q),f\rangle=\displaystyle\frac{(aq;q)_{\infty}}{2(abq^{2};q)_{\infty}}
\displaystyle\sum_{k=0}^{\infty}
  \displaystyle\frac{(bq;q)_{k}}
  {(q;q)_{k}}(aq)^{k}\langle\delta_{q^{k}},f\rangle
 +\frac{K}{2}  \int_{0}^{q^{-1}}x^{\frac{\ln a}{\ln
q}}\displaystyle\frac{(qx;q)_{\infty}}{(bqx;q)_{\infty}} f(x)dx $,\\
 $f\in \mathcal{P}$, $0<q<1$, $0<a <q^{-1}$, $b\in]-\infty,1]\setminus\{0\} $, where
$K^{-1}=\displaystyle\int_{0}^{q^{-1}}x^{\frac{\ln a}{\ln
q}}\displaystyle\frac{(qx;q)_{\infty}}{(bqx;q)_{\infty}}dx$, \\
\tsep{1mm}  $\langle
\mathcal{U}(a,b,q),f\rangle=\displaystyle\frac{(a^{-1}q^{-1};q^{-1})_{\infty}}
{2(a^{-1}b^{-1}q^{-2};q^{-1})_{\infty}}
\displaystyle\sum_{k=0}^{\infty}
  \displaystyle\frac{(b^{-1}q^{-1};q^{-1})_{k}}
  {(q^{-1};q^{-1})_{k}}(aq)^{-k}\langle\delta_{b^{-1}q^{-k-1}},f\rangle $\\
\tsep{1mm}  $\displaystyle\phantom{\langle \mathcal{U}(a,b,q),f\rangle=}{} +\frac{K}{2} \int_{0}^{b^{-1}}x^{\frac{\ln a}{\ln
q}}\displaystyle\frac{(bx;q^{-1})_{\infty}}{(x;q^{-1})_{\infty}}
f(x)dx $,\\
$f\in \mathcal{P}$, $q>1$, $a>q^{-1}$, $b\geq1$  where
 $K^{-1}=\displaystyle\int_{0}^{b^{-1}}x^{\frac{\ln a}{\ln
q}}\displaystyle\frac{(bx;q^{-1})_{\infty}}{(x;q^{-1})_{\infty}}dx$\bsep{3mm}\\
 \hline  \tsep{1mm} \bsep{1mm} case 1.7 \hspace{1cm} $q$-Charlier-II-${\mathcal{U}(\mu,q)}$ \big($\mu\neq 0$, $\mu\neq q^{-n}$,
$n\geq0$\big)\\
\hline
\tsep{4mm}
$\displaystyle \widehat{\beta}_{n}=\frac{1-(1+q)q^{n}+\mu q^{2n}}
{(1-\mu q^{2n-1})(1-\mu q^{2n+1})} q^{n-1}$, $ n\geq 0,$ \\
\tsep{1mm} $\displaystyle \widehat{\gamma}_{n+1}=-q^{3n}\frac{(1-q^{n+1})(1-\mu q^{n})}
{(1-\mu q^{2n}) (1-\mu q^{2n+1})^{2}(1-\mu q^{2n+2})}$, $n\geq 0,$ \\
\tsep{1mm} $H_{q}(x(x-\mu^{-1}q^{-1})\mathcal{U}(\mu,q))-(\mu
q(q-1))^{-1}\{(\mu q-1)x-1\} \mathcal{U}(\mu,q)=0$,
\\
\tsep{1mm} $(\mathcal{U}(\mu,q))_{n}=(-1)^{n}\displaystyle\frac{q^{\frac{1}{2}n(n-1)}}{(\mu
q;q)_{n}}$, $ n \geq 0$, \\
$\displaystyle \langle \mathcal{U}(\mu,q),f\rangle=\frac{1}
{(\mu^{-1}q^{-1};q^{-1})_{\infty}} \sum_{k=0}^{\infty}
  \displaystyle\frac{q^{-\frac{1}{2}k(k+1)}}
  {(q^{-1};q^{-1})_{k}}(-\mu^{-1})^{k}\langle\delta_{\mu^{-1}q^{-k-1}},f\rangle$,
$ f\in \mathcal{P}$, $q>1$, $ \mu < 0 $\\
\hline
\end{tabular}
\end{table}

\begin{table}[t]\centering\small
\centerline{\it Continuation of Table~{\rm \ref{table1}}.}\vspace{1mm}
 \begin{tabular}{|p{15.5cm}|}
 \hline  \tsep{1mm}\bsep{0.5mm} case 1.8 \hspace{1 cm}  Generalized Stieltjes--Wigert  ${\mathcal{S}(\omega,q)}$   \big($\omega\neq q^{-n}$, $n\geq0$\big)
\\
\hline \tsep{2mm}
$\widehat{\beta}_{n}=\{(1+q)q^{-n}-q-\omega\}q^{-n-\frac{3}{2}}$, $n\geq 0$, \\
\tsep{1mm} $\widehat{\gamma}_{n+1}=(1-q^{n+1})(1-\omega q^{n})q^{-4n-4}$, $n\geq 0$, \\
\tsep{1mm} $H_{q}(x(x+\omega q^{-\frac{3}{2}})\mathcal{S}(\omega,q))-(
q-1)^{-1}\big\{x+(\omega-1)q^{-\frac{3}{2}}\big\}\mathcal{S}(\omega,q)=0$,
\\
\tsep{1mm} $(\mathcal{S}(\omega,q))_{n}=q^{-\frac{1}{2}n(n+2)}(\omega;q)_{n}$,
$n \geq 0$,\\
\tsep{1mm} $\mathcal{S}(\omega,q)= (\omega^{-1};q^{-1})_{\infty}
\displaystyle\sum_{k=0}^{\infty}
  \displaystyle\frac{\omega^{-k}}
  {(q^{-1};q^{-1})_{k}}\delta_{-\omega q^{-k-\frac{3}{2}}}$,
$q>1$, $\omega >1$,\\
\tsep{1mm} $\langle
\mathcal{S}(\omega,q),f\rangle=K\displaystyle\int_{0}^{\infty}
\frac{x^{\frac{\ln \omega}{\ln
q}-1}}{(-q^{\frac{3}{2}}\omega^{-1}x;q)_{\infty}} f(x) dx$,  $f\in
\mathcal{P}$, $0<q<1$, $0<\omega <1$,\\
\tsep{1mm}  where  $K^{-1}=\displaystyle\int_{0}^{\infty} \frac{x^{\frac{\ln
\omega}{\ln
q}-1}}{(-q^{\frac{3}{2}}\omega^{-1}x;q)_{\infty}}dx$  is  given by  \eqref{2.10}, \\
\tsep{1mm} $\langle\mathcal{S}(\omega,q),f\rangle=K_{\omega}\displaystyle\int_{q^{-\frac{1}{2}}|\omega|}^{\infty}
(-q^{-\frac {1}{2}}|\omega|x^{-1};q)_{\infty}
\exp\left(-\frac{\ln^{2}x}{2\ln q^{-1}}\right)f(x) dx$,\\
\tsep{1mm} $\phantom{\langle\mathcal{S}(\omega,q),f\rangle=}{}f\in
\mathcal{P}$,  $0<q<1$, $\omega \leq 0$,\\
 $\phantom{\langle\mathcal{S}(\omega,q),f\rangle=}{}$where $K_{\omega}^{-1}=\displaystyle\int_{q^{-\frac{1}{2}}|\omega|}^{\infty}
(-q^{-\frac {1}{2}}|\omega|x^{-1};q)_{\infty}
\exp\left(-\frac{\ln^{2}x}{2\ln q^{-1}}\right)dx$, \\
\tsep{1mm}  $\phantom{\langle\mathcal{S}(\omega,q),f\rangle=}{}$in   particular  $\displaystyle K_{0}=\sqrt{\frac{q}{2\pi \ln q^{-1}}}$ \bsep{2mm}\\
\hline
\end{tabular}\vspace{-2mm}
\end{table}

\begin{table}[t]\centering\small
\caption{Limiting cases.}\label{table2} \vspace{1mm}

\begin{tabular}{|p{15.5cm}|}
  \hline
  \tsep{1mm}\bsep{0.5mm} $H_{q}$-classical linear form \\
   \hline
\tsep{1mm}\bsep{1mm}  case 2.1 \hspace{1cm} $q$-analogue  of  Laguerre  ${\mathbf{L}(\alpha,q)}$  ($\alpha\neq
-[n]_{q}-1$, $n\geq0$)\\
\hline
\tsep{1mm} $\widehat{\beta}_{n}=q^{n}\bigl\{(1+q^{-1})[n]_{q}+1+\alpha\bigr\}$, $n\geq 0,$ \\
\tsep{1mm} $\widehat{\gamma}_{n+1}=q^{2n}[n+1]_{q}\bigl\{[n]_{q}+1+\alpha\bigr\}$, $n\geq 0,$ \\
\tsep{1mm} $H_{q}(x\mathbf{L}
   (\alpha,q))+(x-1-\alpha)\mathbf{L}(\alpha,q)=0$\bsep{1mm}\\
\hline\tsep{1mm}\bsep{1mm} case 2.2 \hspace*{1cm} $q$-analogue  of
  Bessel ${\mathbf{B}(\alpha,q)}$ ($\alpha\neq \frac{1}{2} (q-1)^{-1}$,   $\alpha\neq
-\frac{1}{2}[n]_{q}$, $n\geq0$)\\
\hline
\tsep{3mm} $\displaystyle \widehat{\beta}_{n}=-2 q^{n}\frac{2\alpha
+(1+q^{-1})[n-1]_{q}-q^{-1}
[2n]_{q}}{(2\alpha+[2n-2]_{q})(2\alpha+[2n]_{q})}$, $n\geq 0$, \\
$\displaystyle \widehat{\gamma}_{n+1}=-4 q^{3n}
\frac{[n+1]_{q}(2\alpha+[n-1]_{q})}
{(2\alpha+[2n-1]_{q})(2\alpha+[2n]_{q})^{2}(2\alpha+[2n+1]_{q})}$, $n\geq 0$, \\
$H_{q}(x^{2}\mathbf{B}(\alpha,q))-2(\alpha
x+1)\mathbf{B}(\alpha,q)=0$\bsep{1mm}\\
\hline   \tsep{2mm}\bsep{2mm} case 2.3 \hspace{1 cm} $q$-analogue  of
  Jacobi  ${\mathbf{J}(\alpha,\beta,q)}$ ($\alpha+\beta\neq
\frac{3-2q}{q-1}$,\\
 $\alpha+\beta\neq -[n]_{q}-2$, $n\geq 0$,  $\beta\neq
-[n]_{q}-1$, $n\geq 0$ et $\alpha+\beta+2-(\beta+1)q^{n}+[n]_{q}\neq
0$, $n\geq 0$)\\
\hline
\tsep{4mm} $\displaystyle \widehat{\beta}_{n}=q^{n-1}\frac{(1+q)(\alpha+\beta+2+[n-1]_{q})
(\beta+1+[n]_{q})-(\beta+1 )(\alpha+\beta+2+[2n]_{q})
}{(\alpha+\beta+2+[2n-2]_{q})(\alpha+\beta+2+[2n]_{q})}$, $n\geq 0,$ \\
\tsep{1mm} $\displaystyle \widehat{\gamma}_{n+1}=q^{2n}\frac{[n+1]_{q}(\alpha+\beta+2+[n-1]_{q})
([n]_{q}+\beta+1)(\alpha+\beta+2-(\beta+1)q^{n}+[n]_{q})}
{(\alpha+\beta+2+[2n-1]_{q})(\alpha+\beta+2+[2n]_{q})^{2}(\alpha+\beta+2+[2n+1]_{q})}$,
$n\geq 0,$ \\
\tsep{1mm} $H_{q}(x(x-1)\mathbf{J}(\alpha,\beta,q))-((\alpha+\beta+2)x-(\beta+1))\mathbf{J}(\alpha,\beta,q)=0$
\bsep{1mm}
\\
 \hline
    \end{tabular}
    \vspace{-1mm}
\end{table}

\section{Preliminary and f\/irst results}
\subsection{Preliminary and notations}
\indent Let $\mathcal{P}$ be the vector space of polynomials with
coef\/f\/icients in $\mathbb{C}$ and let $\mathcal{P}^{\prime }$ be its
topological dual. We denote by $\langle u,f\rangle $ the ef\/fect of $
u\in \mathcal{P}^{\prime }$ on $f\in \mathcal{P}$. In particular, we
denote by $(u)_{n}:=\langle u,x^{n}\rangle$, $n\geq 0$ the moments of $u$. Moreover,
a form (linear functional) u is called symmetric if $(u)_{2n+1}=0$,
$n\geq0.$

Let us introduce some useful operations in
$\mathcal{P}^{\prime }$. For any form $u$, any polynomial $g$ and
any $(a,b,c)\in(\mathbb{C}\setminus \{0\})\times\mathbb{C}^{2}$, we
let $H_{q}u$,  $gu$,  $h_{a}u$,  $\tau_{b}u$, $(x-c)^{-1}u$ and
$\delta_{c}$, be the forms def\/ined by duality
\begin{gather*}
\langle H_{q}u,f \rangle :=-\langle
u,H_{q}f\rangle,\qquad \langle gu,f\rangle:=\langle
u,gf\rangle, \qquad \langle h_{a}u,f\rangle:=\langle
u,h_{a}f\rangle,\qquad f\in \mathcal{P},
\\
\langle\tau_{b}u,f\rangle:=\langle
u,\tau_{-b}f\rangle,\qquad \langle (x-c)^{-1}u,f\rangle
:=\langle u,\theta _{c}f\rangle,\qquad \langle
\delta_{c},f\rangle:=f(c),\qquad f\in \mathcal{P},
\end{gather*}
where  $(H_{q}f)(x)= \frac{f(qx)-f(x)}{(q-1)x}$, $q\in
\widetilde{\mathbb{C}}:=\big\{z\in \mathbb{C}, \; z\neq
0, z^{n}\neq 1, n \geq 1\big\}$ \cite{16,18}, $(h_{a}f)(x)=f(ax)$,
$(\tau_{-b}f)(x)=f(x+b)$, $(\theta_{c}f)(x) = \frac{f(x)
-f(c)}{x-c}$ \cite{24} and it's easy to see that \cite{20,26}
\begin{gather}
H_{q}(fu)=(h_{q^{-1}}f)H_{q}u+q^{-1}(H_{q^{-1}}f)u,\qquad
f\in \mathcal{P}, \qquad u\in \mathcal{P}',\label{2.1}
\\
(x-c)((x-c) ^{-1}u) =u ,%\label{2.2}
\qquad
(x-c) ^{-1} (( x-c)u) =u-(u) _{0}\delta _{c}.\nonumber%\label{2.3}
\end{gather}

Now, we introduce the operator $\sigma:\mathcal{P}\longrightarrow
\mathcal{P}$ def\/ined by $(\sigma f)(x):=f(x^{2})$ for all
$f\in\mathcal{P}$. Consequently, we def\/ine $\sigma u$ by duality
\cite{8,25}
\begin{gather*}
\langle \sigma u,f\rangle:=\langle u,\sigma f\rangle,\qquad
f\in \mathcal{P},\qquad u\in \mathcal{P'}.%\label{2.4}
\end{gather*}
We have the well known formula~\cite{25}
\begin{gather}
f(x) \sigma u=\sigma\big(f\big(x^{2}\big) u\big).\label{2.5}
\end{gather}

 Let $\{B_{n}\}_{n\geq 0}$ be a sequence of monic polynomials
with $\deg B_{n}=n$, $n\geq 0$, the form
$u$ is called \textit{regular} if we can associate with it a
sequence of polynomials $\{B_{n}\}_{n\geq 0}$ such that $ \langle
u,B_{m}B_{n}\rangle =r_{n}\delta _{n,m}$, $n, m\geq
0$; $r_{n}\neq 0$, $n\geq 0.$
The sequence $\{ B_{n}\}_{n\geq 0}$ is then said orthogonal with
respect to~$u$. $\{B_{n}\}_{n\geq 0}$ is an (OPS) and it can be
supposed (MOPS). The sequence $\{B_{n}\}_{n\geq 0}$ fulf\/ills the
recurrence relation
\begin{gather}
 B_{0}(x) =1,\qquad  B_{1}(x) = x-\beta _{0}, \nonumber\\
 B_{n+2}(x) =(x-\beta _{n+1}) B_{n+1}(x) -\gamma
_{n+1}B_{n}(x),\qquad \gamma _{n+1}\neq 0,\qquad n\geq 0.
\label{2.6}
\end{gather}
When $u$ is regular, $\{B_{n}\}_{n\geq 0}$ is a symmetrical (MOPS)
if and only if $\beta_{n}=0$, $ n\geq0$.

Lastly, let us recall the following standard expressions \cite{8,11,20}
\begin{gather*}
(a;q)_{0}:=1,\qquad
(a;q)_{n}:=\prod_{k=1}^{n}\big(1-aq^{k-1}\big), \qquad n\geq 1,% \label{2.7}
\\
(a;q)_{\infty}:=\prod_{k=0}^{\infty}\big(1-aq^{k}\big), \qquad |q|<1,%\label{2.8}
\end{gather*}
the $q$-binomial theorem \cite{4,17}
\begin{gather}
\sum_{k=0}^{\infty}\frac{(a;q)_{k}}{(q;q)_{k}}
z^{k}=\frac{(az;q)_{\infty}}{(z;q)_{\infty}},\qquad |z|<1, \qquad |q|<1, \nonumber\\ %\label{2.9}\\
 \int_{0}^{\infty}
t^{x-1} \frac{(-at;q)_{\infty}}{(-t;q)_{\infty}} dt\nonumber\\
\qquad{} =\left\{
     \begin{array}{ll}
       \displaystyle\frac{\pi}{\sin(\pi x)}
 \frac{(a;q)_{\infty}}{(aq^{-x};q)_{\infty}}
 \frac{(q^{1-x};q)_{\infty}}{(q;q)_{\infty}},\quad &
x\in
\mathbb{R}_{+}\setminus\mathbb{N}, \ |a|< q^{x}, \  0<q<1,  \vspace{1mm}\\
\displaystyle\frac{(-q)^{m}}{1-q^{m}}
 \frac{(q^{-1};q^{-1})_{m}}{(aq^{-1};q^{-1})_{m}}
\ln(q^{-1}),\quad & x=m \in \mathbb{N}^{\star},\ |a|<q^{m}, \ 0<q<1.
     \end{array}
   \right.\label{2.10}
\end{gather}

\vspace{-1mm}

\subsection[Some results about the $H_{q}$-semiclassical character]{Some results about the $\boldsymbol{H_{q}}$-semiclassical character}

  A form $u$ is called $H_{q}$-semiclassical when it is
regular and there exist two polynomials $\Phi$ and $\Psi $, $\Phi$
monic, $\deg \Phi=t\geq 0$, $\deg \Psi=p\geq 1$ such that{\samepage
\begin{gather}
H_{q}(\Phi u) +\Psi u=0, \label{2.11}
\end{gather}
the corresponding orthogonal polynomial sequence $\{B_{n}\}_{n\geq
0}$ is
called $H_{q}$-semiclassical \cite{21}.}

The $H_{q}$-semiclassical character is kept by a dilation~\cite{21}. In
fact, let $\{a^{-n}(h_{a}B_{n})\}_{n\geq 0}$, $a\neq 0$; when $u$
satisf\/ies~\eqref{2.11}, then $h_{a^{-1}}u$ fulf\/ills the $q$-analog of the
distributional equation of Pearson type
\begin{gather*}
H_{q}\big(a^{-t}\Phi(ax)h_{a^{-1}}u\big) +a^{1-t}\Psi (ax)
h_{a^{-1}}u=0, %\label{2.12}
\end{gather*}
and the recurrence coef\/f\/icients of~\eqref{2.6} are
\begin{gather*}
\frac{\beta_{n}}{a} ,\qquad
\frac{\gamma_{n+1}}{a^{2}},\qquad n\geq 0.%\label{2.13}
\end{gather*}
Also, the $H_{q}$-semiclassical form $u$ is said to be of class
$s=\max(p-1,t-2)\geq 0$ if and only if~\cite{21}
\begin{gather}
\prod_{c\in \mathcal{Z} _{\Phi }}\bigl\{\bigl| q (h_{q}\Psi)(c)
+(H_{q}\Phi)(c) \bigr| +\bigl| \langle u,q (\theta _{cq}\Psi)  +
(\theta _{cq}\circ\theta _{c} \Phi) \rangle \bigl| \bigr\}>0, \label{2.14}
\end{gather}
where $\mathcal{Z}_{\Phi }$ is the set of zeros of $\Phi$. In
particular, when $s=0$ the form $u$ is usually called
\textit{$H_{q}$-classical} (Al-Salam--Carlitz, big $q$-Laguerre,
$q$-Meixner, Wall, \dots)~\cite{20}.

\begin{lemma}[\cite{21}]\label{lemma1} Let $u$ be a symmetrical
$H_{q}$-semiclassical form of class $s$ satisfying \eqref{2.11}.
The following statements holds
\begin{enumerate}\itemsep=0pt
\item[$i)$] If $s$ is odd then the polynomial $\Phi$ is odd and $\Psi$ is
even.
\item[$ii)$] If $s$ is even then the polynomial $\Phi$ is even and $\Psi$ is
odd.
\end{enumerate}
\end{lemma}
In the sequel we are going to use some $H_{q}$-classical forms~\cite{20},
resumed in Table~\ref{table1} (canonical cases: 1.1--1.8) and Table~\ref{table2} (limiting
cases: 2.1--2.3). %, in the appendix of this contribution.
In fact, when
$q\rightarrow 1$ in results of Table~\ref{table2}, we recover the classical
Laguerre $\mathcal{L}(\alpha)$,   Bessel $\mathcal{B}(\alpha)$ and
$h_{-\frac{1}{2}}\circ \tau_{-1}\mathcal{J}(\alpha,\beta)$
respectively where $\mathcal{J}(\alpha,\beta)$ is the Jacobi
classical form~\cite{24}.

Moreover in what follows we are going to use the logarithmic function denoted by $\textrm{Log}:
\mathbb{C}\setminus \{0\} \longrightarrow \mathbb{C}$ def\/ined by
\begin{gather*}
\textrm{Log}\,z=\ln |z|+i \,\textrm{Arg}\,z,\qquad z \in
\mathbb{C}\setminus \{0\} ,\qquad  -\pi <
\textrm{Arg}\, z\leq \pi,
\end{gather*}
$\textrm{Log}$ is the principal branch of $\log$ and includes $\ln:
\mathbb{R}^{+}\setminus \{0\} \longrightarrow \mathbb{R}$ as a
special case. Consequently, the principal branch of the square root
is
\begin{gather*}
\sqrt{z}=\sqrt{|z|}\,
\textrm{e}^{i\,\frac{\textrm{Arg\,z}}{2}},\qquad z \in
\mathbb{C}\setminus \{0\},\qquad  -\pi <
\textrm{Arg}\, z\leq \pi.
\end{gather*}

\subsection{On quadratic decomposition of a symmetrical regular form}

Let $u$ be a symmetrical regular form and $\{B_{n}\}_{n\geq 0}$ be
its MOPS satisfying \eqref{2.6} with $\beta_{n}=0$, $n\geq 0.$ It is very
well known (see \cite{8,25}) that
\begin{gather*}
B_{2n}(x)=P_{n}\big(x^{2}\big)
,\qquad B_{2n+1}(x)=xR_{n}\big(x^{2}\big),
\qquad n\geq 0,%\label{2.15}
\end{gather*}
where $\{P_{n}\}_{n\geq 0}$ and $\{R_{n}\}_{n\geq 0}$ are the two
MOPS related to the regular form $\sigma u$ and $x \sigma u$
respectively. In fact, \cite{8,25}
\begin{gather*}
\text{$u$ is regular $\Leftrightarrow$ $\sigma u$ and $x \sigma u$ are regular},\\
\text{$u$ is positive def\/inite $\Leftrightarrow$ $\sigma u$ and $x \sigma
u$ are positive def\/inite}.
\end{gather*}
Furthermore, taking
\begin{gather*}
P_{0}(x) =1,\qquad P_{1}(x) =
x-\beta _{0}^{P},\nonumber \\
P_{n+2}(x) =\big(x-\beta _{n+1}^{P}\big) P_{n+1}(x) -\gamma _{n+1}^{P}P_{n}(x),\qquad \gamma _{n+1}^{P}\neq
0,\qquad n\geq 0,
%\label{2.16}
\end{gather*}
and
\begin{gather*}
R_{0}(x) =1,\qquad R_{1}(x) = x-\beta _{0}^{R},\nonumber \\
R_{n+2}(x) =\big(x-\beta _{n+1}^{R}\big) R_{n+1}(x) -\gamma _{n+1}^{R}R_{n}(x),\qquad \gamma _{n+1}^{R}\neq
0,\qquad n\geq 0,%\label{2.17}
\end{gather*}
we get \cite{8,25}
\begin{gather}
\beta_{0}^{P}=\gamma_{1},\nonumber \\
\beta_{n+1}^{P}=\gamma_{2n+2}+\gamma_{2n+3}\hspace{0.1cm},\qquad n\geq
0,\nonumber\\
\gamma_{n+1}^{P} =\gamma_{2n+1} \gamma_{2n+2},\qquad n\geq 0,
\label{2.18}
\end{gather}
and
\begin{gather}
\beta_{n}^{R}=\gamma_{2n+1}+\gamma_{2n+2},\qquad n\geq 0,\nonumber\\
\gamma_{n+1}^{R} =\gamma_{2n+2} \gamma_{2n+3},\qquad n\geq 0.
\label{2.19}
\end{gather}
Consequently,
\begin{gather}
\gamma_{1}=\beta_{0}^{P},\qquad
\gamma_{2}= \frac{\gamma_{1}^{P}}{\beta_{0}^{P}},\nonumber \\
\gamma_{2n+1}=\beta_{0}^{P} \frac{\prod\limits_{k=1}^{n}\gamma_{k}^{R}}
{\prod\limits_{k=1}^{n}\gamma_{k}^{P}},\qquad
\gamma_{2n+2}=\frac{1}{\beta_{0}^{P}}
\frac{\prod\limits_{k=1}^{n+1}\gamma_{k}^{P}}
{\prod\limits_{k=1}^{n}\gamma_{k}^{R}}
,\qquad n\geq 1.
\label{2.20}
\end{gather}

\begin{proposition}\label{proposition1} Let $u$ be a symmetrical regular form.
\begin{enumerate}\itemsep=0pt
\item[$(i)$] The moments of $u$ are
\begin{gather}
(u)_{2n}=(\sigma
u)_{n} ,\qquad (u)_{2n+1}=0,
\qquad n\geq 0.\label{2.21}
\end{gather}

\item[$(ii)$] If $\sigma u$ has the discrete representation
\begin{gather}
\sigma u=\sum_{k=0}^{\infty}\rho_{k}\delta_{\tau_{k}},
\qquad \sum_{k=0}^{\infty}\rho_{k}=1, \label{2.22}
\end{gather}
then a possible discrete measure of $u$ is
\begin{gather}
u=\sum_{k=0}^{\infty}\rho_{k}
\frac{\delta_{\sqrt{\tau_{k}}}+\delta_{_{-}\sqrt{\tau_{k}}}}{2}.\label{2.23}
\end{gather}

\item[$(iii)$] If $u$ is positive definite and $\sigma u$ has the integral
representation
\begin{gather}
\langle \sigma u,f\rangle=\int_{0}^{\infty}V(x) f(x) dx,
\qquad f\in \mathcal{P},
\qquad \int_{0}^{\infty}V(x)dx=1, \label{2.24}
\end{gather}
then, a possible integral representation of $u$ is
\begin{gather}
\langle u,f\rangle=\int_{-\infty}^{\infty}|x|V(x^{2}) f(x) dx,
\qquad f\in \mathcal{P}.\label{2.25}
\end{gather}
\end{enumerate}
\end{proposition}

\begin{proof} $(i)$ is a consequence from the def\/inition
of the quadratic operator $\sigma$.

For $(ii)$ taking into account \eqref{2.21}, \eqref{2.22} we get
\begin{gather*}
(u)_{2n}=(\sigma u)_{n}=\sum_{k=0}^{\infty}\rho_{k}
(\sqrt{\tau_{k}})^{2n}=\sum_{k=0}^{\infty}\rho_{k}
\frac{(\sqrt{\tau_{k}})^{2n}+(_{-}\sqrt{\tau_{k}})^{2n}}{2}.
\end{gather*}
But
\begin{gather*}
(u)_{2n+1}=0=\sum_{k=0}^{\infty}\rho_{k}
\frac{(\sqrt{\tau_{k}})^{2n+1}+(_{-}\sqrt{\tau_{k}})^{2n+1}}{2}.
\end{gather*}
Hence the desired result \eqref{2.23} holds.

For $(iii)$ consider $f\in\mathcal{P}$ and let us split up the
polynomial $f$ accordingly to its even and odd parts
\begin{gather}
f(x)=f^{\rm e}\big(x^{2}\big)+x f^{\rm o}\big(x^{2}\big).\label{2.26}
\end{gather}
Therefore since $u$ is a symmetrical form
\begin{gather}
\langle u,f(x)\rangle=\langle u,f^{\rm e}\big(x^{2}\big)\rangle=\langle \sigma
u,f^{\rm e}(x)\rangle. \label{2.27}
\end{gather}
From \eqref{2.26} we get
\begin{gather}
f^{\rm e}(x)=\frac{f(\sqrt{x})+f(_{-}\sqrt{x})}{2},\qquad x\in
\mathbb{R}_{+}.\label{2.28}
\end{gather}
By \eqref{2.24} and according to \eqref{2.27}, \eqref{2.28} we recover the
representation in \eqref{2.25}.
\end{proof}

\section[Symmetrical $H_{\sqrt{q}}$-semiclassical orthogonal polynomials of class one]{Symmetrical $\boldsymbol{H_{\sqrt{q}}}$-semiclassical orthogonal polynomials\\ of class one}

\begin{lemma} \label{lemma2} We have
\begin{gather}
\sigma(H_{q}u)=(q+1) H_{q^{2}}(\sigma(xu)) ,\qquad u \in
\mathcal{P}'.\label{3.1}
\end{gather}
\end{lemma}

\begin{proof}  From the def\/inition of $H_{q}$ we get
\begin{gather*}
(H_{q}(\sigma f))(x)=(q+1)x (\sigma (H_{q^{2}}f))(x) ,\qquad
f\in \mathcal{P}.
\end{gather*}
Therefore, $ \forall \,f \in\mathcal{P}$,
\begin{gather*}
  \langle \sigma(H_{q}u),f\rangle  = \langle
H_{q}u,\sigma f\rangle  =-\langle u,(q+1)x
\sigma(H_{q^{2}}f)\rangle \\
\phantom{\langle \sigma(H_{q}u),f\rangle}{}
=-\langle (q+1)\sigma (xu),H_{q^{2}}f\rangle  =\langle (q+1) H_{q^{2}}(\sigma (xu)),f\rangle.
\end{gather*}
Thus the desired result.
\end{proof}

\begin{lemma}\label{lemma3} Let $u$ be a
symmetrical $H_{\sqrt{q}}$-semiclassical form of class one. There
exist two polynomials $\varphi$ and $\psi$, $\varphi$ monic, with
$\deg \varphi \leq 1$ and $\deg \psi=1$, such that
\begin{gather}
H_{\sqrt{q}}\big(x \varphi\big(x^{2}\big)u\big)+\psi\big(x^{2}\big) u=0.\label{3.2}
\end{gather}
\end{lemma}

\begin{proof} The result is a consequence from the def\/inition of the
class and Lemma~\ref{lemma1}.
\end{proof}

\begin{corollary}\label{corollary1} Let $u$ be a
symmetrical $H_{\sqrt{q}}$-semiclassical form of class one
satisfying \eqref{3.2}; then $\sigma u$ et $x \sigma u$ are
$H_{q}$-classical satisfying respectively the following $q$-analog
of the distributional equation of Pearson type
\begin{gather}
H_{q}(x\varphi(x) \sigma u)+\frac{1}{\sqrt{q}+1}\psi(x)
\sigma u=0,\label{3.3}\\
H_{q}(x\varphi(x) (x\sigma u))+q^{-1}\left(\frac{1}{\sqrt{q}+1}\psi(x)-\varphi(x)\right)
(x\sigma u)=0.\label{3.4}
\end{gather}
\end{corollary}

\begin{proof} First, $\sigma u$ and $x \sigma u$  are regular because $u$ is symmetrical and regular.
Applying the quadratic operator $\sigma$ to \eqref{3.2} and taking into
account \eqref{3.1} we get
\begin{gather*}
(\sqrt{q}+1)
H_{q}\bigl(\sigma\bigl(x^{2}\varphi\big(x^{2}\big)u\bigr)\bigr)+\sigma\bigl(\psi\big(x^{2}\big)u\bigl)=0.
\end{gather*}
By \eqref{2.5} we get \eqref{3.3}.
Now, multiplying both sides of \eqref{3.3} by $q^{-1}x$, using the
identity in \eqref{2.1}, this yields to~\eqref{3.4}.
\end{proof}

Regarding Table~\ref{table1} (cases 1.1--1.8), Table~\ref{table2} (cases 2.1--2.3) and the
$q$-analog of the distributional equation of Pearson type
\eqref{3.3}, \eqref{3.4}, we consider the following situations for the polynomial
$\varphi$ in order to get a $H_{\sqrt{q}}$-semiclassical form from a
$H_{q}$-classical
\begin{alignat*}{3}
& {\bf A.} \ \varphi(x)=1 \ (\textrm{cases} \ 1.1, 1.2, 1.3, 1.4, 2.1);\qquad &&
{\bf B.} \ \varphi(x)=x \ (\textrm{cases} \ 1.5, 2.2); & \\
& {\bf C.} \ \varphi(x)=x-1 \ (\textrm{case} \ 2.3);\qquad && {\bf D.} \ \varphi(x)=x-b^{-1}q^{-1} \ (\textrm{case} \ 1.6); & \\
& {\bf E.} \ \varphi(x)=x-\mu^{-1}q^{-1} \ (\textrm{case} \ 1.7); \qquad && {\bf F.} \ \varphi(x)=x+\omega
q^{-\frac{3}{2}} \ (\textrm{case} \ 1.8).&
\end{alignat*}

\textbf{A.} In the case \textbf{$\varphi(x)=1$} the $q$-analog of
the distributional equation of Pearson type \eqref{3.3}, \eqref{3.4} are
\begin{gather}
H_{q} (x \sigma u )+\frac{1}{\sqrt{q}+1}\psi(x) \sigma
u=0,\label{3.5}\\
H_{q} (x (x\sigma u))+q^{-1}\left(\frac{1}{\sqrt{q}+1}\psi(x)-1\right) (x\sigma
u)=0.\label{3.6}
\end{gather}

\textbf{A$_{1}$.} If {\bf
$\psi(x)=(\sqrt{q}+1)(x-1-\alpha)$} the $q$-analogue of
the Laguerre form $\mathbf{L}(\alpha,q)$,  $\alpha\neq
-[n]_{q}-1$, $n\geq0$ (case~2.1 in Table~\ref{table2}) satisfying
\begin{gather*}
H_{q}(x\mathbf{L}(\alpha,q))+(x-1-\alpha)\mathbf{L}(\alpha,q)=0.%\label{3.7}
\end{gather*}
Comparing with \eqref{3.5}, \eqref{3.6} we get
\begin{gather}
\sigma u=\mathbf{L}(\alpha,q),\qquad \alpha\neq -[n]_{q}-1,\qquad n\geq0,\label{3.8}
\end{gather}
and
\begin{gather}
x\sigma u=(1+\alpha) \mathbf{L}\big(q^{-1}(\alpha+2)-1,q\big),\qquad
\alpha\neq -[n]_{q}-1,\qquad n\geq0.\label{3.9}
\end{gather}
Taking into account the recurrence coef\/f\/icients (see case~2.1 in
Table~\ref{table2}), by virtue of \eqref{3.8}, \eqref{3.9} and \eqref{2.18}, \eqref{2.19} we get for $n
\geq 0$
\begin{gather*}
    \beta_{n}^{P}=q^{n}\big\{\big(1+q^{-1}\big)[n]_{q}+1+\alpha\big\}, \nonumber \\
    \gamma_{n+1}^{P}=q^{2n}[n+1]_{q} \{[n]_{q}+1+\alpha\}, \nonumber \\
    \beta_{n}^{R}= q^{n+1}\big\{\big(1+q^{-1}\big)[n]_{q}+q^{-1}(2+\alpha)\big\},\nonumber \\
    \gamma_{n+1}^{R}= q^{2n+2}[n+1]_{q} \big\{[n]_{q}+q^{-1}(2+\alpha)\big\}.
  %\label{3.10}
\end{gather*}
With the relation $[k-1]_{q}=q^{-1}[k]_{q}-q^{-1}$, $k\geq 1$ the
system \eqref{2.20} becomes for $n \geq 0$
\begin{gather}
    \gamma_{2n+1}=q^{n}([n]_{q}+1+\alpha), \qquad
    \gamma_{2n+2}=q^{n} [n+1]_{q}.
  \label{3.11}
\end{gather}
Writing $\alpha=\mu-\frac{1}{2}$, $\mu \neq
-[n]_{q}-\frac{1}{2}$, $n\geq0$ and denoting the symmetrical form
$u$ by $\mathcal{H}(\mu,q)$ we get the following result:

\begin{proposition} \label{proposition2} The symmetrical form $\mathcal{H}(\mu,q)$ satisfies the following properties:
\begin{enumerate}\itemsep=0pt
\item[$1)$] The recurrence coefficient $\gamma_{n+1}$ satisfies
\eqref{3.11}.

\item[$2)$] $\mathcal{H}(\mu,q)$ is regular if and only if
$\mu \neq
-[n]_{q}-\frac{1}{2}$, $n\geq0$.

\item[$3)$] $\mathcal{H}(\mu,q)$ is positive definite if and only
if
$q>0$, $\mu>-\frac{1}{2}$.

\item[$4)$] $\mathcal{H}(\mu,q)$ is a $H_{\sqrt{q}}$-semiclassical
form of class one for $\mu\neq
\frac{1}{\sqrt{q}(\sqrt{q}+1)}-\frac{1}{2}$, $\mu \neq
-[n]_{q}-\frac{1}{2}$, $n\geq0$ satisfying the $q$-analog of the
distributional equation of Pearson type
\begin{gather}
H_{\sqrt{q}} (x\mathcal{H}(\mu,q))+(\sqrt{q}+1)\left(
x^{2}-\mu-\frac{1}{2}\right)\mathcal{H}(\mu,q)=0.\label{3.12}
\end{gather}
\end{enumerate}
\end{proposition}

\begin{proof} The results in 1), 2) and 3) are straightforward from
\eqref{3.11}. For 4), it is clear that $\mathcal{H}(\mu,q)$ satisf\/ies~\eqref{3.12}; in
this case and by virtue of~\eqref{2.14}, we are going to prove that the
class of $\mathcal{H}(\mu,q)$ is exactly one for $\mu\neq
\frac{1}{\sqrt{q}(\sqrt{q}+1)}-\frac{1}{2}$, $\mu \neq
-[n]_{q}-\frac{1}{2}$, $n\geq0$.
Denoting $\Phi(x)=x$, $\Psi(x)=(\sqrt{q}+1)\bigl(
x^{2}-\mu-\frac{1}{2}\bigr)$, we have accordingly to \eqref{2.14}, on one
hand
\begin{gather*}
\sqrt{q}\bigl(h_{\sqrt{q}}\Psi\bigr)(0)+\bigl(H_{\sqrt{q}}\Phi\bigr)(0)
=1-\sqrt{q}(\sqrt{q}+1)\left(\mu+\frac{1}{2}\right)\neq0,
\end{gather*}
and on the other hand by $(\theta_{0}\Psi)(x)=(\sqrt{q}+1)x$ and
$(\theta_{0}^{2}\Phi)(x)=0$,
\begin{gather*}
\langle \mathcal{H}(\mu,q), \sqrt{q}
\theta_{0}\Psi+\theta_{0}^{2}\Phi\rangle=0,
\end{gather*}
taking into account that $u$ is a symmetrical form.
\end{proof}

\begin{remark} The symmetrical form $\mathcal{H}(\mu,q)$, $\mu\neq
\frac{1}{\sqrt{q}(\sqrt{q}+1)}-\frac{1}{2}$, $\mu \neq
-[n]_{q}-\frac{1}{2}$,  $n\geq0$ is the $q$-analogue of the
generalized Hermite one~\cite{12} (when $q\rightarrow 1$ we recover the
generalized Hermite form $\mathcal{H}(\mu)$ (see~\eqref{1.3}) which is a
symmetrical semiclassical form of class one for   $\mu\neq0$,
$\mu \neq -n-\frac{1}{2}$, $n\geq0$ \cite{1,8,15,26}).
\end{remark}

\textbf{A$_{2}$.} If {\bf $\psi(x)=-(\sqrt{q}-1)^{-1}(x+1)$} the
form $\mathcal{U}$ that satisf\/ies the $q$-analog of the
distributional equation of Pearson type (see case 1.1 in Table~\ref{table1})
\begin{gather*}
H_{q}(x\mathcal{U})-(q-1)^{-1}(x+1)\mathcal{U}=0.%\label{3.13}
\end{gather*}
Comparing with \eqref{3.5}, \eqref{3.6} we get
\begin{gather}
\sigma u=\mathcal{U}, \label{3.14}
\end{gather}
and
\begin{gather}
x\sigma u=-h_{q}\mathcal{U}.\label{3.15}
\end{gather}
Taking into account \eqref{3.14}, \eqref{3.15}, \eqref{2.18}, \eqref{2.19} and the case 1.1 in
Table~\ref{table1} we obtain for $n \geq 0$
\begin{gather*}
    \beta_{n}^{P}=\big\{1-(1+q)q^{n}\big\}q^{n-1}, \qquad
    \gamma_{n+1}^{P}=\big(q^{n+1}-1\big)q^{3n}, \nonumber \\
    \beta_{n}^{R}=\big\{1-(1+q)q^{n}\big\}q^{n}, \qquad
    \gamma_{n+1}^{R}=\big(q^{n+1}-1\big)q^{3n+2}.
 % \label{3.16}
\end{gather*}
Consequently, the system \eqref{2.20} becomes for $n \geq 0$
\begin{gather}
    \gamma_{2n+1}=-q^{2n}, \qquad
    \gamma_{2n+2}=\big(1-q^{n+1}\big)q^{n}.
  \label{3.17}
\end{gather}

\begin{proposition}\label{proposition3} The symmetrical form $u$ satisfies the following properties:
\begin{enumerate}\itemsep=0pt
\item[$1)$] The recurrence coefficient $\gamma_{n+1}$ satisfies
\eqref{3.17}.
\item[$2)$] $u$ is regular for any $q\in \widetilde{\mathbb{C}}$.
\item[$3)$] $u$ is a $H_{\sqrt{q}}$-semiclassical form of class
one satisfying
\begin{gather}
H_{\sqrt{q}}(x u)-(\sqrt{q}-1)^{-1}\big(x^{2}+1\big)u=0.\label{3.18}
\end{gather}

\item[$4)$] The moments of $u$ are
\begin{gather*}
(u)_{2n}=(-1)^{n} q^{\frac{1}{2}n(n-1)},\qquad (u)_{2n+1}=0   , \qquad n\geq 0.%\label{3.19}
\end{gather*}

\item[$5)$] we have the following discrete representation
\begin{gather*}
u=\sum_{k=0}^{\infty} \frac{(-1)^{k}
q^{-k^{2}}s(k)}{(q^{-1};q^{-1})_{k}}
\frac{\delta_{iq^{\frac{k}{2}}}+\delta_{-iq^{\frac{k}{2}}}}{2}, \qquad q>1.%\label{3.20}
\end{gather*}
\end{enumerate}
\end{proposition}

\begin{proof} The results in 1), 2) are obvious from \eqref{3.17}.
For 3), it is clear that $u$ satisf\/ies~\eqref{3.18}. Denoting $\Phi(x)=x$,
$\Psi(x)=-(\sqrt{q}-1)^{-1}\big(x^{2}+1\big)$, we have \eqref{2.14}
\begin{gather*}
\sqrt{q}\bigl(h_{\sqrt{q}}\Psi\bigr)(0)+\bigl(H_{\sqrt{q}}\Phi\bigr)(0)=\frac{1}{1-\sqrt{q}}
\neq 0,
\qquad
\langle u, \sqrt{q} \theta_{0}\Psi+\theta_{0}^{2}\Phi\rangle=0.
\end{gather*}
Therefore, $u$ is of class one.
The results in 4) and 5) are consequence from \eqref{2.21}--\eqref{2.23} and
those for $\mathcal{U}$ (case 1.1 in Table~\ref{table1}).
\end{proof}

\textbf{A$_{3}$.} If {\bf
$\psi(x)=-(aq)^{-1}(\sqrt{q}-1)^{-1} (x-1+aq )$} the little
$q$-Laguerre form $\mathcal{L}(a,q)$, $a\neq 0$, $a\neq
q^{-n-1}$, $n\geq0$ (case~1.2 in Table~\ref{table1}) satisfying
\begin{gather*}
H_{q}(x\mathcal{L}(a,q))-(aq)^{-1}(q-1)^{-1}(x-1+aq)\mathcal{L}(a,q)=0.%\label{3.21}
\end{gather*}
With \eqref{3.5}, \eqref{3.6} we obtain
\begin{gather}
\sigma u=\mathcal{L}(a,q),\qquad  a\neq 0,\qquad  a\neq q^{-n-1},\qquad
n\geq0, \label{3.22}
\end{gather}
and
\begin{gather}
x\sigma u=(1-aq)\mathcal{L}(aq,q),\qquad a\neq 0,\qquad a\neq
q^{-n-1},\qquad n\geq0.\label{3.23}
\end{gather}
By virtue of the recurrence coef\/f\/icients of little $q$-Laguerre
polynomials in Table~\ref{table1}, case~1.2, the relations in~\eqref{3.22}, \eqref{3.23} and
\eqref{2.18}, \eqref{2.19} we get for $n \geq 0$
\begin{gather*}
    \beta_{n}^{P}=\big\{1+a-a(1+q)q^{n}\big\}q^{n}, \nonumber  \\
    \gamma_{n+1}^{P}=a\big(1-q^{n+1}\big)\big(1-aq^{n+1}\big)q^{2n+1}, \nonumber \\
    \beta_{n}^{R}=\big\{1+aq-a(1+q)q^{n+1}\big\}q^{n}, \nonumber \\
    \gamma_{n+1}^{R}=a\big(1-q^{n+1}\big)\big(1-aq^{n+2}\big)q^{2n+2}.
 % \label{3.24}
\end{gather*}
Therefore \eqref{2.20} becomes for $n \geq 0$
\begin{gather}
    \gamma_{2n+1}=q^{n}\big(1-aq^{n+1}\big), \qquad
    \gamma_{2n+2}=a q^{n+1}\big(1-q^{n+1}\big).
  \label{3.25}
\end{gather}
Comparing with~\cite{2}, $u$ is a symmetrical case of the Al-Salam--Verma
form, $u:=\mathcal{SV}(a,q)$. From~\eqref{3.25}, it is easy to see that
$\mathcal{SV}(a,q)$ is regular if and only if   $a\neq 0$, $a\neq q^{-n-1}$, $n\geq0$. Also, $\mathcal{SV}(a,q)$ is positive
def\/inite if and only if  $0<q<1$, $0<a<q^{-1}$  or
 $q>1$, $a<0$.

\begin{proposition}\label{proposition4} The form $\mathcal{SV}(a,q)$ is a $H_{\sqrt{q}}$-semiclassical form of class one for
$a\neq 0$, $a\neq q^{-\frac{1}{2}}$, $a\neq q^{-n-1}$, $n\geq0$ satisfying
\begin{gather}
H_{\sqrt{q}} (x\mathcal{SV}(a,q))
-(aq)^{-1}(\sqrt{q}-1)^{-1}\bigl(x^{2}-1+aq\bigr)\mathcal{SV}(a,q)=0.\label{3.26}
\end{gather}
The moments are
\begin{gather}
(\mathcal{SV}(a,q))_{2n}=(aq;q)_{n},\qquad (\mathcal{SV}(a,q))_{2n+1}=0  ,
\qquad n\geq 0,\label{3.27}
\end{gather}
and the orthogonality relation can be represented
\begin{gather}
 \langle\mathcal{SV}(a,q),f\rangle=\frac{(aq;q)_{\infty}}{2}
\displaystyle\sum_{k=0}^{\infty}\frac{(aq)^{k}}{(q;q)_{k}}\biggl\langle
\frac{\delta_{q^{\frac{k}{2}}}+\delta_{-q^{\frac{k}{2}}}}{2},f\biggl\rangle \nonumber
\\
{}+\frac{K}{2}\int_{-q^{-\frac{1}{2}}}^{q^{-\frac{1}{2}}}|x|^{2\frac{\ln
a}{\ln q}+1} (qx^{2};q)_{\infty} f(x) dx , \qquad f \in \mathcal{P}, \qquad
0<q<1, \qquad 0<a<q^{-1},\label{3.28}
\end{gather}
with
\begin{gather*}
K^{-1}=q^{-\frac{\ln a}{\ln
q}-1} \int_{0}^{1}x^{\frac{\ln a}{\ln q}} (x;q)_{\infty}
dx,
\end{gather*}
 and
\begin{gather}
\mathcal{SV}(a,q)=\frac{1}{(a;q^{-1})_{\infty}} \sum_{k=0}^{\infty}
\frac{q^{-\frac{1}{2}k(k-1)}(-a)^{k}}{(q^{-1};q^{-1})_{k}}
\frac{\delta_{q^{\frac{k}{2}}}+\delta_{-q^{\frac{k}{2}}}}{2}, \qquad
q>1, \qquad a<0.\label{3.29}
\end{gather}
\end{proposition}
\begin{proof} It is direct that the form $\mathcal{SV}(a,q)$ satisf\/ies
the $q$-analog of the distributional equation of Pearson type~\eqref{3.26}. Denoting $\Phi(x)=x$, $\Psi(x)=-(aq)^{-1}(\sqrt{q}-1)^{-1}\bigl(x^{2}-1+aq\bigr)$, we have~\eqref{2.14}
\begin{gather*}
\sqrt{q}\bigl(h_{\sqrt{q}}\Psi\bigr)(0)+\bigl(H_{\sqrt{q}}\Phi\bigr)(0)=
\frac{a^{-1}q^{-\frac{1}{2}}-1}{\sqrt{q}-1}\neq 0,
\qquad
\langle \mathcal{SV}(a,q), \sqrt{q}
\theta_{0}\Psi+\theta_{0}^{2}\Phi\rangle=0,
\end{gather*}
from which we get that $\mathcal{SV}(a,q)$ is of class one because
$a\neq
0$, $ a\neq q^{-\frac{1}{2}}$, $a\neq q^{-n-1}$, $n\geq0$.
The results mentioned in \eqref{3.27}--\eqref{3.29} are easily obtained from
those well known the properties of the little $q$-Laguerre from
(case 1.2 in Table~\ref{table1}) and \eqref{2.21}--\eqref{2.25}.
\end{proof}

\begin{remark} The regular form $\mathcal{SV}(q^{-\frac{1}{2}},q)$ is
the discrete $\sqrt{q}$-Hermite form which is
$H_{\sqrt{q}}$-classical~\cite{20}.
\end{remark}

\textbf{A$_{4}$.} If {\bf
$\psi(x)=-b^{-1}(\sqrt{q}-1)^{-1}\bigl(q^{-1}x+b-1\bigr)$} the Wall
form $\mathcal{W}(b,q)$,  $b\neq 0$, $b\neq q^{-n}$, $n\geq0$
(case~1.3 in Table~\ref{table1}) that satisf\/ies
\begin{gather*}
H_{q}(x\mathcal{W}(b,q))-b^{-1}(q-1)^{-1}(q^{-1}x+b-1)\mathcal{W}(b,q)=0.%\label{3.30}
\end{gather*}
In accordance of \eqref{3.5}, \eqref{3.6} we get
\begin{gather*}
\sigma u=\mathcal{W}(b,q), \qquad b\neq 0,\qquad b\neq q^{-n},\qquad
n\geq0,%\label{3.31}
\end{gather*}
and
\begin{gather*}
x\sigma u=q(1-b)\mathcal{W}(bq,q), \qquad b\neq 0,\qquad b\neq
q^{-n},\qquad n\geq0.%\label{3.32}
\end{gather*}
We recognize the Brenke type symmetrical regular form
$\mathcal{Y}(b,q)$ \cite{8,9,10}. In~\cite{13} it is proved that
$\mathcal{Y}(b,q)$ is $H_{\sqrt{q}}$-semiclassical of class one for
$b\neq 0$, $b\neq \sqrt{q}$, $b\neq q^{-n}$, $n\geq0$
satisfying
\begin{gather}
H_{\sqrt{q}} (x\mathcal{Y}(b,q))-b^{-1}\big(q^{\frac{1}{2}}-1\big)^{-1}\big\{q^{-1}
x^{2}+b-1\big\}\mathcal{Y}(b,q)=0.\label{3.33}
\end{gather}
Also in that work, moments, discrete and integral representations
are established.

\begin{remark} Likewise, from \eqref{3.33} it is easy to see that
$h_{\frac{1}{\sqrt{q}}}\mathcal{Y}(\sqrt{q},q)$ is the
$H_{\sqrt{q}}$-classical discrete $\sqrt{q}$-Hermite  form~\cite{20}.
\end{remark}

\textbf{A$_{5}$.} If {\bf
$\psi(x)=(\sqrt{q}-1)^{-1}q^{-\alpha-1}\bigl(x+b-q^{\alpha+1}\bigr)$}
the generalized $q^{-1}$-Laguerre $\mathcal{U}^{(\alpha)}(b,q)$
form,  $b\neq 0$, $b\neq q^{n+1+\alpha}$, $n\geq 0 $
and its $q$-analog of the distributional equation of Pearson type
(case~1.4 in Table~\ref{table1})
\begin{gather*}
H_{q}(x\mathcal{U}^{(\alpha)}(b,q))+(q-1)^{-1}q^{-\alpha-1}(x+b-q^{\alpha+1})
\mathcal{U}^{(\alpha)}(b,q)=0.%\label{3.34}
\end{gather*}
By \eqref{3.5}, \eqref{3.6} we deduce the following relationships
\begin{gather}
\sigma u=\mathcal{U}^{(\alpha)}(b,q), \qquad b\neq 0,\qquad b\neq
q^{n+1+\alpha},\qquad n\geq0   , \label{3.35}
\\
x\sigma u=\big(q^{\alpha+1}-b\big)\mathcal{U}^{(\alpha+1)}(b,q), \qquad b\neq
0,\qquad b\neq q^{n+1+\alpha},\qquad n\geq0 . \label{3.36}
\end{gather}
From Table~\ref{table1}, case~1.4, the relations in \eqref{3.35}, \eqref{3.36} and
\eqref{2.18}, \eqref{2.19} we get for $n \geq 0$
\begin{gather*}
    \beta_{n}^{P}=\big\{1-q^{-n-1}+q^{-1}\big(1-bq^{-n-\alpha}\big)\big\}q^{2n+\alpha+1}, \nonumber\\
    \gamma_{n+1}^{P}=\big(1-q^{-n-1}\big)\big(1-bq^{-n-1-\alpha}\big)q^{4n+2\alpha+3}, \nonumber \\
    \beta_{n}^{R}=\big\{1-q^{-n-1}+q^{-1}\big(1-bq^{-n-\alpha-1}\big)\big\}q^{2n+\alpha+2}, \nonumber \\
    \gamma_{n+1}^{R}=\big(1-q^{-n-1}\big)\big(1-bq^{-n-2-\alpha}\big)q^{4n+2\alpha+5}.  %\label{3.37}
\end{gather*}
Thus, for $n \geq 0$
\begin{gather*}
    \gamma_{2n+1}=\big(1-bq^{-n-1-\alpha}\big)q^{2n+\alpha+1}, \qquad
    \gamma_{2n+2}=\big(1-q^{-n-1}\big)q^{2n+\alpha+2}.
 % \label{3.38}
\end{gather*}
Consequently, the symmetrical form $u:=u(\alpha,b,q)$ is regular if
and only if $b \neq 0$, $b\neq q^{n+1+\alpha}$, $n\geq0$. It
is positive def\/inite for $\alpha\in \mathbb{R}$, $q>1$, $b<q^{\alpha+1}$.

\begin{proposition} \label{proposition5} The symmetrical form $u$ is a $H_{\sqrt{q}}$-semiclassical form of class one
for $b\neq 0$, $b\neq q^{n+1+\alpha}$, $n\geq0$, $\alpha\in \mathbb{R}$ satisfying
\begin{gather*}
H_{\sqrt{q}} (x
u)+q^{-\alpha-1}\big(q^{\frac{1}{2}}-1\big)^{-1}\big\{
x^{2}+b-q^{\alpha+1}\big\}u=0.%\label{3.39}
\end{gather*}
Moreover, we have the following identities
\begin{gather}
(u)_{2n}=(-b)^{n}\big(b^{-1}q^{\alpha+1};q\big)_{n}
 ,\qquad (u)_{2n+1}=0  , \qquad n\geq 0,\label{3.40}\\
\langle
u,f\rangle=K\int_{-\infty}^{\infty}\frac{|x|^{2\alpha-2\frac{\ln
b}{\ln q}+1}} {(-b^{-1}x^{2};q^{-1})_{\infty}} f(x) dx, \label{3.41}
\end{gather}
for $f \in \mathcal{P}$,  $\alpha\in \mathbb{R}$, $q>1$,
$0<b<q^{\alpha+1}$, with
\begin{gather*}
K^{-1}= \int_{0}^{\infty}\frac{x^{\alpha-\frac{\ln
b}{\ln q}}}
{(-b^{-1}x;q^{-1})_{\infty}} dx
\end{gather*}
 is given by \eqref{2.10},
\begin{gather}
u=\frac{1}{(b^{-1}q^{\alpha};q^{-1})_{\infty}} \sum_{k=0}^{\infty}
\frac{q^{-\frac{1}{2}k(k-1)}}{(q^{-1};q^{-1})_{k}}
(-b^{-1}q^{\alpha})^{k}
\frac{\delta_{\sqrt{-b}q^{\frac{k}{2}}}+\delta_{-\sqrt{-b}q^{\frac{k}{2}}}}{2} ,  \label{3.42}
\end{gather}
for $\alpha\in \mathbb{R}$, $q>1$, $b<0$, and
\begin{gather}
u=\big(b^{-1}q^{\alpha+1};q\big)_{\infty} \sum_{k=0}^{\infty}
\frac{(b^{-1}q^{\alpha+1})^{k}}{(q;q)_{k}}
\frac{\delta_{i\sqrt{b}q^{\frac{k}{2}}}+\delta_{-i\sqrt{b}q^{\frac{k}{2}}}}{2},  \label{3.43}
\end{gather}
for $\alpha\in \mathbb{R}$, $0<q<1$, $b>q^{\alpha+1}$.
\end{proposition}

\begin{proof} First, let us obtain the class of the form; denoting
\begin{gather*}
\Phi(x)=x, \qquad \Psi(x)=(\sqrt{q}-1)^{-1}q^{-\alpha-1}\bigl(x^{2}+b-q^{\alpha+1}\bigr),
\end{gather*}
we have
\begin{gather*}
\sqrt{q}\bigl(h_{\sqrt{q}}\Psi\bigr)(0)+\bigl(H_{\sqrt{q}}\Phi\bigr)(0)=\frac{b
q^{-\alpha-\frac{1}{2}}-1}{\sqrt{q}-1} \neq 0,
\qquad \langle u, \sqrt{q} \theta_{0}\Psi+\theta_{0}^{2}\Phi\rangle=0,
\end{gather*}
for $b\neq 0$, $b\neq q^{n+1+\alpha}$, $n\geq0$, $\alpha\in \mathbb{R}$. Thus, $u$ is of class one.
The identities given in \eqref{3.40}--\eqref{3.43} are easily obtained from the
properties of the generalized $q^{-1}$-Laguerre
$\mathcal{U}^{(\alpha)}(b,q)$ form (Table~\ref{table1}, case~1.4) and
\eqref{2.21}--\eqref{2.25}.
\end{proof}

\textbf{B.} In the case $\varphi(x)=x$ the
$q$-analog of the distributional equation of Pearson type
\eqref{3.3}, \eqref{3.4} are
\begin{gather}
H_{q}\bigl(x^{2} \sigma u\bigr)+\frac{1}{\sqrt{q}+1}\psi(x) \sigma
u=0,\label{3.44}\\
H_{q}\big(x^{2} (x\sigma
u)\big)+q^{-1}\left\{\frac{1}{\sqrt{q}+1}\psi(x)-x\right\} (x\sigma
u)=0.\label{3.45}
\end{gather}

\textbf{B$_{1}$.} If {\bf $\psi(x)=-2(\sqrt{q}+1) (\alpha x+1)$} the
$q$-analogue of the Bessel form (case 2.2 in Table~\ref{table2}), the form
$\mathbf{B}(\alpha,q)$, $\alpha\neq \frac{1}{2}(q-1)^{-1}$,
$\alpha\neq -\frac{1}{2}[n]_{q}$, $n\geq0$ satisfying
\begin{gather*}
H_{q}\big(x^{2}\mathbf{B}(\alpha,q)\big)-2(\alpha
x+1)\mathbf{B}(\alpha,q)=0.%\label{3.46}
\end{gather*}
Thus, comparing with \eqref{3.44}, \eqref{3.45}, we get
\begin{gather*}
\sigma u=\mathbf{B}(\alpha,q), \qquad \alpha\neq
\frac{1}{2}(q-1)^{-1} ,\qquad \alpha\neq -\frac{1}{2}[n]_{q},\qquad
n\geq0,%\label{3.47}
\end{gather*}
and
\begin{gather*}
x\sigma u=-\alpha^{-1}
h_{q^{-1}}\mathbf{B}(q^{-1}(\alpha+\frac{1}{2}),q),\,\, \alpha\neq
\frac{1}{2}(q-1)^{-1},\qquad \alpha\neq -\frac{1}{2}[n]_{q},\qquad
n\geq0.%\label{3.48}
\end{gather*}
By the recurrence coef\/f\/icients in case 2.2 of Table~\ref{table2}, the relations in
\eqref{3.44}, \eqref{3.45} and \eqref{2.18}, \eqref{2.19} we get for $n \geq 0$
\begin{gather*}
    \beta_{n}^{P}=-2q^{n}
    \frac{2\alpha + (1+q^{-1})[n-1]_{q}-q^{-1}[2n]_{q}}
    {(2\alpha+[2n-2]_{q})(2\alpha+[2n]_{q})}  , \nonumber \\
    \gamma_{n+1}^{P}=-4q^{3n}
    \frac{[n+1]_{q}(2\alpha+[n-1]_{q})}
    {(2\alpha+[2n-1]_{q})(2\alpha+[2n]_{q})^{2}(2\alpha+[2n+1]_{q})}  , \nonumber \\
      \beta_{n}^{R}=-2q^{n-1}
    \frac{(2\alpha+1)q^{-1} + (1+q^{-1})[n-1]_{q}-q^{-1}[2n]_{q}}
    {((2\alpha+1)q^{-1}+[2n-2]_{q})((2\alpha+1)q^{-1}+[2n]_{q})}  , \nonumber \\
    \gamma_{n+1}^{R}=-4q^{3n-2}
    \frac{[n+1]_{q}((2\alpha+1)q^{-1}+[n-1]_{q})}
    {((2\alpha+1)q^{-1}\!+[2n-1]_{q})
    ((2\alpha+1)q^{-1}\!+[2n]_{q})^{2}((2\alpha+1)q^{-1}\!+[2n+1]_{q})} .  %\label{3.49}
\end{gather*}
By the relation $[k-1]_{q}=q^{-1}[k]_{q}-q^{-1}$, $k\geq 1$,
\eqref{2.20} leads to for $n\geq 0$
\begin{gather}
  \gamma_{1}=- \frac{1}{\alpha}  , \qquad
  \gamma_{2n+2}=2q^{2n}  \frac{[n+1]_{q}}   {(2\alpha+[2n]_{q})(2\alpha+[2n+1]_{q})} , \nonumber \\
  \gamma_{2n+3}=-2q^{n+1}  \frac{(2\alpha+[n]_{q})}
    {(2\alpha+[2n+1]_{q})(2\alpha+[2n+2]_{q})} .
  \label{3.50}
\end{gather}
We put $\alpha= \frac{\nu+1}{2}$, $\nu \neq
\frac{2-q}{q-1}$, $\nu \neq -[n]_{q}-1$, $n\geq0$ and denote the
symmetrical form $u$ by $\mathcal{B}[\nu,q]$. From~\eqref{3.50} the form
$\mathcal{B}[\nu,q]$ is regular if and only if  $\nu \neq
\frac{2-q}{q-1}$, $\nu \neq -[n]_{q}-1$, $n\geq0$. Also,
it is quite straightforward to deduce that the symmetrical form
$\mathcal{B}[\nu,q]$ is $H_{\sqrt{q}}$-semiclassical of class one
for $\nu \neq \frac{2-q}{q-1}$, $\nu \neq -[n]_{q}-1$, $n\geq0$
satisfying the $q$-analog of the distributional equation of Pearson
type
\begin{gather*}
H_{\sqrt{q}} \bigl(x^{3}\mathcal{B}[\nu,q]\bigr)-2(\sqrt{q}+1)\left(
\frac{\nu+1}{2} x^{2}+1\right)\mathcal{B}[\nu,q]=0.%\label{3.51}
\end{gather*}
\begin{remark} The symmetrical form
$h_{(2\sqrt{2})^{-1}}\mathcal{B}[\nu,q]$, $\nu \neq \frac{2-q}{q-1}
$, $\nu \neq -[n]_{q}-1$, $n\geq0$ is the $q$-analogue of
the symmetrical form $\mathcal{B}[\nu]$~\cite{14} (when $q\rightarrow 1$
we recover the symmetrical semiclassical~$\mathcal{B}[\nu]$, $\nu
\neq -n-1$, $n\geq0$ of class one, see~\eqref{1.5}). Also, for any
parameter $\alpha\neq -n-1$, $n\geq0$ the symmetrical form
$h_{(2 \sqrt{1+\sqrt{q}})^{-1}}
\mathcal{B}[-\frac{q^{-\alpha-1}-1}{q-1}-1,q]$ appears in~\cite{33}.
\end{remark}

\textbf{B$_{2}$.} If {\bf $\psi(x)=-(aq)^{-1}(\sqrt{q}-1)^{-1}
((1+aq) x-1)$ } the Alternative $q$-Charlier
$\mathcal{A}(a,q)$ form with $a\neq 0$, $a\neq -q^{-n}$, $n\geq 0 $ that satisf\/ies (case~1.5 in Table~\ref{table1})
\begin{gather*}
H_{q}\big(x^{2}\mathcal{A}(a,q)\big)-(aq)^{-1}(q-1)^{-1} \big((1+aq)
x-1\big) \mathcal{A}(a,q)=0.%\label{3.52}
\end{gather*}
Thus
\begin{gather*}
\sigma u=\mathcal{A}(a,q), \qquad a\neq 0,\qquad a\neq -q^{-n},\qquad
n\geq0 , %\label{3.53}
\end{gather*}
and
\begin{gather*}
x\sigma u=\frac{1}{1+aq}\mathcal{A}(aq,q), \qquad a\neq 0,\qquad a\neq -q^{-n},\qquad n\geq0 . %\label{3.54}
\end{gather*}
The systems \eqref{2.18}, \eqref{2.19} are for $n \geq 0$
\begin{gather*}
    \beta_{n}^{P}=q^{n}  \frac{1+aq^{n-1}+aq^{n}-aq^{2n}}
    {(1+aq^{2n-1})(1+aq^{2n+1})}  ,  \nonumber \\
    \gamma_{n+1}^{P}=a q^{3n+1} \frac{(1-q^{n+1})(1+a q^{n})}
    {(1+a q^{2n})(1+a q^{2n+1})^{2}(1+a q^{2n+2})}  , \nonumber \\
    \beta_{n}^{R}= q^{n}  \frac{1+aq^{n}+aq^{n+1}-aq^{2n+1}}
    {(1+aq^{2n})(1+aq^{2n+2})}  ,\nonumber \\
    \gamma_{n+1}^{R}= a q^{3n+2} \frac{(1-q^{n+1})(1+a q^{n+1})}
    {(1+a q^{2n+1})(1+a q^{2n+2})^{2}(1+a q^{2n+3})} ,%\label{3.55}
\end{gather*}
from which we get for $n\geq 0$
\begin{gather*}
    \gamma_{2n+1}=q^{n} \frac{1+aq^{n}}
    {(1+aq^{2n})(1+aq^{2n+1})}  , \qquad
    \gamma_{2n+2}= a q^{2n+1} \frac{1-q^{n+1}}
    {(1+aq^{2n+1})(1+aq^{2n+2})}  .%\label{3.56}
\end{gather*}
Consequently, the symmetrical form $u=u(a,q)$ is regular if and only
if  $a \neq 0$, $a\neq -q^{-n}$, $n\geq0$. It is positive
def\/inite for $0<q<1$, $a>0$. Also, $u$ is
$H_{\sqrt{q}}$-semiclassical of class one for $a\neq 0$, $a\neq
-q^{-n}$, $n\geq0$ satisfying the $q$-analog of the distributional
equation of Pearson type
\begin{gather*}
H_{\sqrt{q}}\big(x^{3} u\big)-(aq)^{-1}(\sqrt{q}-1)^{-1}\big(
(1+aq)x^{2}-1\big)u=0.%\label{3.57}
\end{gather*}
After some straightforward computations, we get the following
representations for the moments and the orthogonality
\begin{gather*}
(u)_{2n}=\frac{1}{(-aq;q)_{n}}
,\qquad (u)_{2n+1}=0  , \qquad n\geq
0,%\label{3.58}
\\
\langle u,f\rangle=q^{\frac{1}{2}(\frac{\ln a}{\ln
q}+\frac{1}{2})^{2}}\frac{(-a^{-1};q)_{\infty}}{ \sqrt{2\pi \ln
q^{-1}}}\displaystyle\int_{-\infty}^{\infty}|x|^{2\frac{\ln a}{\ln
q}}(qx^{2};q)_{\infty}\exp\left(-2\frac{\ln^{2}|x|}{\ln q^{-1}}\right)
f(x)dx,  %\label{3.59}
\end{gather*}
for $f \in \mathcal{P}$, $0<q<1$, $a>0$, and
\begin{gather*}
u=\frac{1}{(-aq;q)_{\infty}} \sum_{k=0}^{\infty}
\frac{a^{k}q^{\frac{1}{2}k(k+1)}}{(q;q)_{k}}
\frac{\delta_{-q^{\frac{k}{2}}}+\delta_{q^{\frac{k}{2}}}}{2}  ,
\qquad 0<q<1, \qquad a>0.%\label{3.60}
\end{gather*}

\textbf{C.} In the case  $\varphi(x)=x-1$ the
$q$-analogue of Jacobi form (case~2.3 in Table~\ref{table2}), therefore the
$q$-analog of the distributional equation of Pearson type
\eqref{3.3}, \eqref{3.4} become
\begin{gather*}
H_{q}(x(x-1) \sigma u)-\big((\alpha+\beta+2)x-(\beta+1)\big)
\sigma u=0 ,%\label{3.61}
\end{gather*}
and
\begin{gather*}
H_{q}\big(x(x-1) (x\sigma
u)\big)-q^{-1}\big((\alpha+\beta+3)x-(\beta+2)\big) (x\sigma
u)=0. %\label{3.62}
\end{gather*}
Consequently,
\begin{gather}
\sigma u= \mathbf{J}(\alpha,\beta,q), \label{3.63}
\\
x\sigma
u=\frac{\beta+1}{\alpha+\beta+2}\mathbf{J}\big(q^{-1}(\alpha+1)-1,q^{-1}(\beta+2)-1,q\big)
\label{3.64}
\end{gather}
with the constraints
\begin{gather}
\alpha+\beta\neq \frac{3-2q}{q-1}, \qquad \alpha+\beta\neq
-[n]_{q}-2, \qquad \beta\neq -[n]_{q}-1,\nonumber\\
\alpha+\beta+2-(\beta+1)q^{n}+[n]_{q}\neq 0,\qquad n\geq 0.\label{3.65}
\end{gather}
By Table~\ref{table2} and \eqref{3.63}, \eqref{3.64}, the systems \eqref{2.18}, \eqref{2.19} give for $n
\geq 0$
\begin{gather*}
    \beta_{n}^{P}=q^{n-1}\frac{(1+q)(\alpha+\beta+2+[n-1]_{q})
(\beta+1+[n]_{q})-(\beta+1 )(\alpha+\beta+2+[2n]_{q})
}{(\alpha+\beta+2+[2n-2]_{q})(\alpha+\beta+2+[2n]_{q})} ,\\ %\label{3.66} \\
    \gamma_{n+1}^{P}= q^{2n}\frac{[n+1]_{q}(\alpha+\beta+2+[n-1]_{q})
([n]_{q}+\beta+1)(\alpha+\beta+2-(\beta+1)q^{n}+[n]_{q})}
{(\alpha+\beta+2+[2n-1]_{q})(\alpha+\beta+2+[2n]_{q})^{2}(\alpha+\beta+2+[2n+1]_{q})},
\nonumber\\
      \beta_{n}^{R}=q^{n-1}\frac{(1+q)(\alpha+\beta+2+[n]_{q})
(\beta+1+[n+1]_{q})-(\beta+2 )(\alpha+\beta+2+[2n+1]_{q})
}{(\alpha+\beta+2+[2n-1]_{q})(\alpha+\beta+2+[2n+1]_{q})}  , \nonumber \\
    \gamma_{n+1}^{R}= q^{2n+1}\frac{[n\!+\!1]_{q}(\alpha\!+\!\beta\!+\!2\!+\![n]_{q})
([n+1]_{q}\!+\!\beta\!+\!1)(\alpha\!+\!\beta\!+\!2\!-\!(\beta+2)q^{n}\!+\![n+1]_{q})}
{(\alpha+\beta+2+[2n]_{q})(\alpha+\beta+2+[2n+1]_{q})^{2}(\alpha+\beta+2+[2n+2]_{q})}  . \nonumber
\end{gather*}
Using the above results and the relations
\begin{gather*}
[k-1]_{q}=q^{-1}[k]_{q}-q^{-1}   , \qquad
[k]_{q}=q^{k-1}+[k-1]_{q}  , \qquad k\geq 1
\end{gather*}
we deduce from \eqref{2.20} for $n\geq 0$
\begin{gather}
  \gamma_{2n+1}=q^{n}  \frac{(\alpha+\beta+2+[n-1]_{q})(\beta+1+[n]_{q})}
  {(\alpha+\beta+2+[2n-1]_{q})(\alpha+\beta+2+[2n]_{q})}  , \nonumber   \\
  \gamma_{2n+2}=q^{n} [n+1]_{q} \frac{\alpha+\beta+2-(\beta+1)q^{n}+[n]_{q}}
    {(\alpha+\beta+2+[2n]_{q})(\alpha+\beta+2+[2n+1]_{q})}   .\label{3.67}
\end{gather}
We denote the symmetrical form $u$ by $\mathcal{G}(\alpha,\beta,q)$.
From~\eqref{3.67} the symmetrical form $\mathcal{G}(\alpha,\beta,q)$ is
regular if and only if the conditions in~\eqref{3.65} hold. It is
$H_{\sqrt{q}}$-semiclassical of class one for $\alpha+\beta\neq
\frac{3-2q}{q-1}$, $\alpha+\beta\neq -[n]_{q}-2$, $\beta\neq -[n]_{q}-1$, $\alpha+\beta+2-(\beta+1)q^{n}+[n]_{q}\neq
0$, $n\geq 0$, $\beta\neq \frac{1}{\sqrt{q}(\sqrt{q}+1)}-1$
satisfying
\begin{gather*}
H_{q}\big(x(x^{2}-1)\mathcal{G}(\alpha,\beta,q)\big)-
(\sqrt{q}+1)
\big((\alpha+\beta+2)x^{2}-(\beta+1)\big)\mathcal{G}(\alpha,\beta,q)=0.%\label{3.68}
\end{gather*}

\begin{remark} The symmetrical form $\mathcal{G}(\alpha,\beta,q)$ is the
$q$-analogue of the symmetrical generalized Gegenbauer
$\mathcal{G}(\alpha,\beta)$ form (see~\eqref{1.4}) which is semiclassical
of class one for
$\alpha\neq-n-1$, $\beta\neq-n-1$,
$\beta\neq-\frac{1}{2}$,
$\alpha+\beta\neq-n-1$, $n\geq0$~\cite{1,6}.
\end{remark}

\textbf{D.} In the case
$\varphi(x)=x-b^{-1}q^{-1}$ the little $q$-Jacobi
$\mathcal{U}(a,b,q)$ form (case~1.6 in Table~\ref{table1}). The $q$-analog of the
distributional equation of Pearson type in \eqref{3.3}, \eqref{3.4} become
\begin{gather*}
H_{q}\big(x\big(x-b^{-1}q^{-1}\big) \sigma u\big)+\big(abq^{2}(q-1)\big)^{-1}\big(\big(1-abq^{2}\big)x+aq-1\big) \sigma u=0
,%\label{3.69}
\\
H_{q}\big(x\big(x-b^{-1}q^{-1}\big) (x\sigma u)\big)+\big(abq^{3}(q-1)\big)^{-1}\big(\big(1-abq^{3}\big)x+aq^{2}-1\big)
(x\sigma u)=0. %\label{3.70}
\end{gather*}
Hence
\begin{gather}
\sigma u= \mathcal{U}(a,b,q), \label{3.71}\\
x\sigma u=\frac{1-aq}{1-abq^{2}}\mathcal{U}(aq,b,q) \label{3.72}
\end{gather}
with the constraints
\begin{gather}
ab\neq 0, \qquad a\neq q^{-n-1} , \qquad b\neq q^{-n-1},\qquad ab\neq
q^{-n},\qquad n\geq 0.\label{3.73}
\end{gather}
By Table~\ref{table1} and \eqref{3.71}, \eqref{3.72}, the systems \eqref{2.18}, \eqref{2.19} lead to for
$n \geq 0$
\begin{gather*}
    \beta_{n}^{P}=q^{n} \frac{(1+a)(1+abq^{2n+1})
-a(1+b)(1+q)q^{n}
}{(1-abq^{2n})(1-abq^{2n+2})}  , \nonumber \\
    \gamma_{n+1}^{P}= aq^{2n+1} \frac{(1-q^{n+1})
(1-aq^{n+1})(1-bq^{n+1})(1-abq^{n+1})}
{(1-abq^{2n+1})(1-abq^{2n+2})^{2}(1-abq^{2n+3})}   ,
\nonumber\\
      \beta_{n}^{R}=q^{n} \frac{(1+aq)(1+abq^{2n+2})
-a(1+b)(1+q)q^{n+1}
}{(1-abq^{2n+1})(1-abq^{2n+3})}   ,\nonumber \\
    \gamma_{n+1}^{R}= aq^{2n+2} \frac{(1-q^{n+1})
(1-aq^{n+2})(1-bq^{n+1})(1-abq^{n+2})}
{(1-abq^{2n+2})(1-abq^{2n+3})^{2}(1-abq^{2n+4})}  .
%\label{3.74}
\end{gather*}
Using the above results and~\eqref{2.20} we get for $n \geq 0$
\begin{gather*}
  \gamma_{2n+1}=q^{n} \frac{
(1-aq^{n+1})(1-abq^{n+1})}
{(1-abq^{2n+1})(1-abq^{2n+2})}  , \qquad
  \gamma_{2n+2}= a q^{n+1} \frac{(1-q^{n+1})
(1-bq^{n+1})} {(1-abq^{2n+2})(1-abq^{2n+3})} .
%\label{3.75}
\end{gather*}
Therefore, the symmetrical form $u=u(a,b,q)$ is regular if and only
if the conditions in \eqref{3.73} are satisf\/ied. Further, the form $u$ is
positive def\/inite for $0<q<1$, $0<a<q^{-1}$, $b<1$, $b\neq0$ or $q>1$, $a>q^{-1}$, $b\geq1$. Moreover, by virtue of \eqref{2.14}, the form $u$ is
$H_{\sqrt{q}}$-semiclassical of class one for $ab\neq 0$, $a\neq
q^{-n-1}$, $b\neq q^{-n-1}$, $ab\neq q^{-n}$, $n\geq 0$, $a\neq q^{-\frac{1}{2}}$
\begin{gather*}
H_{\sqrt{q}} \big(x\big(x^{2}-b^{-1}q^{-1}\big)
u\big)+\big(abq^{2}(\sqrt{q}-1)\big)^{-1}\big(\big(1-abq^{2}\big)x^{2}+aq-1\big)u=0.%\label{3.76}
\end{gather*}
Proposition~\ref{proposition1} and the well known representations of the little
$q$-Jacobi form (Table~\ref{table1}) allow us to establish the following
results
\begin{gather*}
(u)_{2n}=\frac{(aq;q)_{n}}{(abq^{2};q)_{n}}
,\qquad (u)_{2n+1}=0  , \qquad  n\geq
0.%\label{3.77}
\end{gather*}
For $f \in \mathcal{P}$, $0<q<1$, $0<a<q^{-1}$, $b<1$, $b\neq 0$,
\begin{gather*}
u= \frac{(aq;q)_{\infty}}{(abq^{2};q)_{\infty}} \sum_{k=0}^{\infty}
\frac{(aq)^{k}(bq;q)_{k})}{(q;q)_{k}}
\frac{\delta_{-q^{\frac{k}{2}}}+\delta_{q^{\frac{k}{2}}}}{2},%\label{3.78}
\end{gather*}
and
\begin{gather*}
\langle u,f\rangle= K \int_{-q^{-\frac{1}{2}}}^{q^{-\frac{1}{2}}}
|x|^{2\frac{\ln a}{\ln q}+1}
\frac{(qx^{2};q)_{\infty}}{(bqx^{2};q)_{\infty}} f(x) dx,
%\label{3.79}
\end{gather*}
with
\begin{gather*}
K^{-1}= \int_{0}^{q^{-1}} x^{\frac{\ln a}{\ln q}}
\frac{(qx;q)_{\infty}}{(bqx;q)_{\infty}}dx.
\end{gather*}
For $f \in \mathcal{P}$, $q>1$, $a>q^{-1}$, $b\geq1$
\begin{gather*}
u=
\frac{(a^{-1}q^{-1};q^{-1})_{\infty}}{(a^{-1}b^{-1}q^{-2};q^{-1})_{\infty}}
\sum_{k=0}^{\infty}
\frac{(aq)^{-k}(b^{-1}q^{-1};q^{-1})_{k})}{(q^{-1};q^{-1})_{k}}
\frac{\delta_{-\sqrt{b^{-1}}q^{-\frac{k+1}{2}}}+
\delta_{\sqrt{b^{-1}}q^{-\frac{k+1}{2}}}}{2},%\label{3.80}
\end{gather*}
and
\begin{gather*}
\langle u,f\rangle= K\int_{-b^{-\frac{1}{2}}}^{b^{-\frac{1}{2}}}
|x|^{2\frac{\ln a}{\ln q}+1}
\frac{(bx^{2};q^{-1})_{\infty}}{(x^{2};q^{-1})_{\infty}} f(x) dx,
%\label{3.81}
\end{gather*}
with
\begin{gather*}
K^{-1}= \int_{0}^{b^{-1}} x^{\frac{\ln a}{\ln q}}
\frac{(bx;q^{-1})_{\infty}}{(x;q^{-1})_{\infty}}dx.
\end{gather*}

\textbf{E.} In the case
$\varphi(x)=x-\mu^{-1}q^{-1}$ the $q$-Charlier-II-form
$\mathcal{U}(\mu,q)$ (case~1.7 in Table~\ref{table1}). From the above assumption
\eqref{3.3}, \eqref{3.4} are
\begin{gather*}
H_{q}\big(x\big(x-\mu^{-1}q^{-1}\big) \sigma u\big)-(\mu q(q-1))^{-1}\big((\mu
q-1)x-1\big) \sigma u=0 ,%\label{3.82}
\\
H_{q}\big(x\big(x-\mu^{-1}q^{-1}\big) (x\sigma u)\big)-(\mu
q(q-1))^{-1}\big(\big(\mu q^{2}-1\big)q^{-1} x-1\big) (x\sigma u)=0.
%\label{3.83}
\end{gather*}
Thus
\begin{gather*}
\sigma u= \mathcal{U}(\mu,q),\qquad \mu\neq 0  , \qquad  \mu\neq
q^{-n}, \qquad n\geq 0,%\label{3.84}
\\
x\sigma u=\frac{1}{\mu q-1}h_{q}\mathcal{U}(\mu q,q),\qquad \mu\neq 0, \qquad \mu\neq q^{-n}, \qquad n\geq 0. %\label{3.85}
\end{gather*}
Consequently, the systems \eqref{2.18}, \eqref{2.19} for $n \geq 0$ are
\begin{gather*}
    \beta_{n}^{P}=q^{n-1} \frac{1
-(1+q)q^{n}+\mu q^{2n}}
{(1-\mu q^{2n-1})(1-\mu q^{2n+1})} , \nonumber \\
    \gamma_{n+1}^{P}= -q^{3n} \frac{(1-q^{n+1})
(1-\mu q^{n})} {(1-\mu q^{2n})(1-\mu q^{2n+1})^{2}(1-\mu q^{2n+2})}
 ,
 \nonumber\\
      \beta_{n}^{R}= q^{n} \frac{1
-(1+q)q^{n}+\mu q^{2n+1}}
{(1-\mu q^{2n})(1-\mu q^{2n+2})}   , \nonumber \\
    \gamma_{n+1}^{R}=-q^{3n+2} \frac{(1-q^{n+1})
(1-\mu q^{n+1})} {(1-\mu q^{2n+1})(1-\mu q^{2n+2})^{2}(1-\mu
q^{2n+3})}   .
%\label{3.86}
\end{gather*}
On account of \eqref{2.20} we have for $n \geq 0$
\begin{gather*}
  \gamma_{2n+1}=-q^{2n} \frac{
1-\mu q^{n}}
{(1-\mu q^{2n})(1-\mu q^{2n+1})}  , \qquad
  \gamma_{2n+2}= q^{n} \frac{1-q^{n+1}}
  {(1-\mu q^{2n+1})(1-\mu q^{2n+2})}  .%\label{3.87}
\end{gather*}
From the last result, the symmetrical form  $u=u(\mu,q)$ is regular
if and only if $\mu\neq 0$, $\mu\neq q^{-n}$, $n\geq 0$.
Moreover, by virtue of \eqref{2.14}, it is clear that $u$ is
$H_{\sqrt{q}}$-semiclassical of class one for $\mu\neq 0$, $\mu\neq q^{-n}$, $n\geq 0$ satisfying
\begin{gather*}
H_{\sqrt{q}} \big(x\big(x^{2}-\mu^{-1}q^{-1}\big) u\big)-(\mu
q(\sqrt{q}-1))^{-1}\big\{(\mu q-1)x^{2}-1\big\}u=0.%\label{3.88}
\end{gather*}
Furthermore, by the same procedure as in {\bf D} we get{\samepage
\begin{gather*}
(u)_{2n}=\frac{(-1)^{n}q^{\frac{1}{2}n(n-1)}}{(\mu q;q)_{n}}
,\qquad (u)_{2n+1}=0  , \qquad n\geq
0,%\label{3.89}
\\
u= \frac{1}{(\mu^{-1}q^{-1};q^{-1})_{\infty}} \sum_{k=0}^{\infty}
\frac{(-\mu^{-1})^{k}q^{-\frac{1}{2}k(k+1)}}{(q^{-1};q^{-1})_{k}}
\frac{\delta_{-i\sqrt{-\mu^{-1}}q^{-\frac{k+1}{2}}}+
\delta_{i\sqrt{-\mu^{-1}}q^{-\frac{k+1}{2}}}}{2}, %\label{3.90}
\end{gather*}
for $q>1$, $\mu<0$.}

\textbf{F.} In the case $\varphi(x)=x+\omega
q^{-\frac{3}{2}}$ the generalized Stieltjes--Wigert form
$\mathcal{S}(\omega,q)$ (case~1.8 in Table~\ref{table1}). From \eqref{3.3}, \eqref{3.4} it
follows
\begin{gather*}
H_{q}\big(x\big(x+\omega q^{-\frac{3}{2}}\big) \sigma u\big)-(q-1)^{-1}\big(x+(\omega-1)q^{-\frac{3}{2}}\big) \sigma u=0,%\label{3.91}
\\
H_{q}\big(x\big(x+\omega q^{-\frac{3}{2}}\big) (x\sigma u)\big)-(q-1)^{-1}\big(x+(\omega q-1)q^{-\frac{5}{2}}\big)
(x\sigma u)=0. %\label{3.92}
\end{gather*}
Thus
\begin{gather*}
\sigma u= \mathcal{S}(\omega,q) , \qquad \omega\neq q^{-n},\qquad
n\geq 0,%\label{3.93}
\\
x\sigma u=(1-\omega)q^{-\frac{3}{2}} h_{q^{-1}}\mathcal{S}(\omega
q,q)  , \qquad \omega\neq q^{-n},\qquad n\geq 0. %\label{3.94}
\end{gather*}
We obtain for $n\geq0$
\begin{gather*}
    \beta_{n}^{P}=\big\{(1+q)q^{-n}-q-\omega\big\}q^{-n-\frac{3}{2}}, \nonumber\\
    \gamma_{n+1}^{P}= \big(1-q^{n+1}\big)\big(1-\omega q^{n}\big)q^{-4n-4},
\nonumber\\
      \beta_{n}^{R}=\big\{(1+q)q^{-n}-q(1+\omega)\big\}q^{-n-\frac{5}{2}}, \nonumber \\
    \gamma_{n+1}^{R}= \big(1-q^{n+1}\big)\big(1-\omega q^{n+1}\big)q^{-4n-6}.
%\label{3.95}
\end{gather*}
Thus, \eqref{2.20} gives for $n\geq0$
\begin{gather}
  \gamma_{2n+1}=q^{-2n-\frac{3}{2}}\big(1-\omega q^{n}\big), \qquad
  \gamma_{2n+2}= q^{-2n-\frac{5}{2}}\big(1-q^{n+1}\big). \label{3.96}
\end{gather}
We recognize the Brenke type symmetrical orthogonal polynomials
\cite{8,9,10}
\begin{gather*}
B_{n}=T_{n}(\cdot ;\omega,q)  , \qquad n \geq 0.
\end{gather*}
We denote $u=\mathcal{T}(w,q)$. Taking into consideration \eqref{3.96},
the symmetrical form $\mathcal{T}(\omega,q)$ is regular if and only
if $\omega\neq q^{-n}$, $n\geq 0$, and it is positive def\/inite for
$0<q<1$, $\omega<1$. Furthermore, it is easy to deduce that
$\mathcal{T}(\omega,q)$ is $H_{\sqrt{q}}$-semiclassical of class one
for $\omega\neq\sqrt{q}$, $\omega\neq q^{-n}$, $n\geq 0$
satisfying the $q$-analog of the distributional equation of Pearson
type
\begin{gather*}
H_{\sqrt{q}} \big(x\big(x^{2}+\omega q^{-\frac{3}{2}}\big)
\mathcal{T}(\omega,q)\big) -(\sqrt{q}-1)^{-1}\big
(x^{2}+(\omega-1)q^{-\frac{3}{2}}\big)
\mathcal{T}(\omega,q)=0.%\label{3.97}
\end{gather*}
Finally, with Proposition~\ref{proposition1} and the properties of the generalized
Stieltjes--Wigert $H_{q}$-classical form (Table~\ref{table1}, case~1.8) we deduce the
following results
\begin{gather*}
(\mathcal{T}(\omega,q))_{2n}=q^{-\frac{1}{2}n(n+2)}(\omega ;q)_{n}
,\qquad (\mathcal{T}(\omega,q))_{2n+1}=0   ,
 n\geq 0,%\label{3.98}
\\
\mathcal{T}(\omega,q)= \big(\omega^{-1};q^{-1}\big)_{\infty}
\sum_{k=0}^{\infty}
\frac{\omega^{-k}}{(q^{-1};q^{-1})_{k}}
\frac{\delta_{-i\sqrt{\omega}q^{-\frac{k}{2}-\frac{3}{4}}}+
\delta_{i\sqrt{\omega}q^{-\frac{k}{2}-\frac{3}{4}}}}{2} , \qquad
q>1  , \qquad \omega>1,%\label{3.99}
\\
\langle\mathcal{T}(\omega,q),f\rangle= K \int_{-\infty}^{\infty}
\frac{|x|^{2\frac{\ln \omega}{\ln
q}-1}}{(-q^{\frac{3}{2}}\omega^{-1}x^{2};q)_{\infty}}f(x) dx, \nonumber\\
f\in \mathcal{P}, \qquad 0<q<1   , \qquad 0<\omega<1, %\label{3.100}
\end{gather*}
with
\begin{gather*}
K^{-1}= \int_{0}^{\infty} \frac{x^{\frac{\ln
\omega}{\ln q}-1}}{(-q^{\frac{3}{2}}\omega^{-1}x;q)_{\infty}} dx
\end{gather*}
 is
given by \eqref{2.10},
\begin{gather*}
\langle \mathcal{T}(0,q),f\rangle=\sqrt{\frac{q}{2\pi \ln q^{-1}}}
\int_{-\infty}^{\infty} |x| \exp\left(\frac{-2 \ln^{2}|x|}{\ln
q^{-1}}\right)f(x) dx   , \qquad 0<q<1.%\label{3.101}
\end{gather*}

\subsection*{Acknowledgments}

The authors are very grateful to the editors and the referees for
the constructive and valuable comments and recommendations.

\pdfbookmark[1]{References}{ref}
\LastPageEnding

\end{document}